\documentclass[12pt]{amsart}

\usepackage{amsmath}
\usepackage{amssymb}
\usepackage{latexsym}
 \newtheorem{thm}{Theorem}[section] 
\newtheorem{cor}[thm]{Corollary} 
\newtheorem{lemma}[thm]{Lemma} 
\newtheorem{prop}[thm]{Proposition} 
\newtheorem{defn}[thm]{Definition} 
\newtheorem{rem}[thm]{Remark} 
\newenvironment{proof*}{\vskip 2mm\noindent {}}{$\Box$ \vskip 2mm}
\numberwithin{equation}{section} 
\newcommand{\C}{{\mathbb{C}}}

\begin{document} 
%\vspace{1cm} 
 \title{Pluricomplex Green and Lempert functions 
for equally weighted
poles} 
                     %\vskip1cm 

                     \author{Pascal J. Thomas \ \ \ Nguyen Van Trao} 
\maketitle

\section{Introduction}

The pluricomplex Green function with several poles introduced by Lelong
\cite{lelong} is one of the most important tools of complex
pluripotential theory. For details we refer the readers to
\cite{Lempert},
\cite{Klimek} and \cite{Coman}.

Let us recall the definition of the pluricomplex Green function with 
several poles. Let $\Omega $ be a domain in $\C^n$, and poles and
weights denoted by
$$
S=\{(a_1,\nu_1) ;...; (a_N,\nu_N)\} \subset \Omega \times \mathbb
R_+,
$$
where $\mathbb R_+=[0,+\infty)$. Define the pluricomplex Green function
\begin{eqnarray*}
\lefteqn{G_S(z):=} \\
&&
\sup \left\{u(z): u\in PSH_-(\Omega), u(x)\le \nu_j \log
\|x-a_j\|+C_j \; \text{when}\; x \to a_j, j=1,...,N \right\}.
\end{eqnarray*}
Note that if $N=1$ we might as well take $\nu_1 =1$,
and the above-mentioned function $g$ is the
pluricomplex Green function with one pole defined for instance in
\cite{Klimek}.

We also recall the definition of Coman's Lempert function
  \cite{Coman}:
$$
\ell_S(z):=\inf \left\{ \sum^N_{j=1}\nu_j\log|\zeta_j|: \exists
\varphi\in \mathcal {O}(\mathbb D,\Omega), \varphi(0)=z,
\varphi(\zeta_j)=a_j, j=1,...,N \right\},
$$
where $\mathbb D$ is the unit disc in $\mathbb C$.

It is easy to see that $\ell_S(z) \ge G_S(z)$ for all $z \in \Omega$.

A remarkable theorem of Lempert \cite{Lempert}
says that equality holds in the case where $\Omega$ is convex and
$N=1$ (and then the weight $\nu_1$ doesn't matter).
Later Coman \cite{Coman} proved with considerable effort
  that this assertion also holds
when $\Omega$ is the unit ball, $N=2$, and the weights are equal. At
the same time he conjectured that the equality might hold for any 
number
of points and
any convex domain in $\mathbb C^n$. Recently, Carlehed and
Wiegerinck \cite{LnotequalG}, \cite{CarlWieg} proved that
Coman's conjecture fails
for the bidisc, with two poles lying on a coordinate axis and distinct
  weights. The main goal of this chapter is to prove that Coman's 
conjecture
does not even hold in the case
when all weights are equal.  The following suggests that this is a more
"natural" case.

Since weights on the Green function are analogous to multiplicities for
zeros,
in view of \cite{LnotequalG}, \cite{CarlWieg}, we focus on "multiple 
poles", that is,
the behavior of the Coman's Lempert function with many poles when
every group of poles tends to some pole. Note that, when we consider 
the
pluricomplex Green function as an elementary solution to the complex
Monge-Amp\`ere operator in several variables, the quantity which we
expect to see being preserved under this limit process is the total
Monge-Amp\`ere mass of the function, which is equal to $\sum_j \nu_j^n$
for a Green function with weights. When a group of $N$ poles, with $N
\neq a^n$ for any integer $a$, clusters to a single point, we cannot
hope to have an usual weighted Green function arise
as limit value for the sequence of Green functions for the separate
poles. Simple examples yield explicit non-isotropic functions, i.e. 
which
are not equivalent to constant multiples of the logarithm of a norm.

Lelong and  Rashkovskii
\cite{lelongandrashkovskii} introduced a generalized pluricomplex Green
function with many poles (see the definition in Section \ref{defns}), 
which allows
for non-isotropic singularities.
We then study the  problem of producing an analogous generalization
of the case $\nu_j=1$ for Coman's
Lempert function.
A motivation is that
we know that Coman's Lempert function is continuous with
respect to $z$ and to its poles when they stay away from each other 
(see
\cite{ThomCCP} for the case of the ball), and we would like to extend 
such
results to singular situations arising from "collisions" of poles.

Unfortunately, this was not fully successful, since our candidate is 
not
in general the limit of the Lempert functions for the natural systems 
of
points which tend to the given "multiple poles". However, we gather
enough information to prove that in some cases involving four points
forming a product set in the bidisc, with all weights equal to one,
equality does not
hold between the Lempert and Green functions.

Along the way, we give partial answers. There is equality between
Lelong and Rashkovskii's Green function and our generalization of 
Coman's
Lempert function in the case of one pole, in the polydisc, with a 
simple enough singularity
(Lemma \ref{L_S=G_S}; some hypothesis about integer multiplicities is 
of course necessary).
We also prove equality between
Lempert and Green functions in the case of
the bidisc in $\mathbb C^2$, when all poles are on a coordinate disk 
and
all multiplicities equal to one ; and also in the natural limit cases 
of
those, when the non-isotropic singularities are all "horizontal",
oriented along the coordinate disc (this is made precise in Theorem
\ref{thm:G=LinthecaseNpoles}). Then, we find the
limit of the Lempert functions in the case when two fixed poles $a_1, 
a_2$ lie on a
coordinate axis, $a_3$ lies on a line orthogonal to this axis at
$a_1$, and $a_3$ tends to $a_1$ (the limit of the corresponding Green
functions is not known in this case).

The organization of the paper is as follows : in Section \ref{defns}, we 
give notations and definitions, 
introduce our generalization of the Lempert function and 
give Lemma \ref{L_S=G_S} as a first motivation of this particular 
definition. In Section \ref{extdiscs}, we generalize to this new Lempert 
functions some of the results of \cite{WikstromAMS}. Section \ref{addone} is 
devoted to the technical details (mostly in one complex variable)
of the proofs of the previous Section. Section \ref{bidisc} provides a few 
positive and negative examples in the bidisc, the latter motivating a 
corrected definition of the generalization of Coman's Lempert function. 
Finally, in
Section \ref{mcex}, using the results of the previous section,
we study the
following situation: the four poles are $(a,0), (b,0), (b,\varepsilon),
(a,\varepsilon)$ and $\varepsilon $ tends to $0$. This
 provides the counterexample to Coman's 
conjecture for single poles (Theorem \ref{NoComan}), 
and also shows that the corrected definition 
still does not yield the limit of the Lempert function under collisions of 
poles (Theorem \ref{limell4<ltilde}), despite a partial positive result 
(Proposition \ref{prop:limsupl4epsilon}). 

\section{Definitions}
\label{defns}

We now introduce some notations.

A set $A \subset \mathbb C^n$ is called $0- circled$ if $x=(x_j)\in A$
implies $x'=(x_je^{i\theta_j})\in A$ for $0\le \theta_j\le 2\pi, 1\le
j\le n.$ We will say that a function $f(x)$ defined on $A,$ is circled
if it is invariant with respect to the rotations $x_j \to x_je^{i
\theta_j}, 1\le j\le n$.

%The notion of $x^0-circled$ sets is defined
%similarly.

In what follows, we will use a special subclass of circled
plurisubharmonic functions \cite{lelongandrashkovskii} $f\in PSH_-(\mathbb D^n)$ that have the
following "conic" property: the convex image $g_f(u)$ of $f$ satisfies
the equation
$$
g_f(Cu)=Cg_f(u), \; \text{for every}\; C>0,
$$
where $x=(x_j)=(\exp(u_j+i\theta_j)); u=(u_j),$ and $g_f(u)=f(x)$. Such
a function $f$ will be called an \emph{indicator}. 

Given a function $f\in PSH(\Omega)$ and a point $x^0 \in \Omega,$
Lelong and Rashkovskii \cite{lelongandrashkovskii} constructed a
function $\Psi_{f,x^0}(y)$ related to local properties of $f$ at $x^0$.
For instance, $\Psi_{f,x^0}(y)\in PSH_-(\mathbb D^n),$ where $ \mathbb
D^n$ being the open unit polydisk in the $\mathbb C^n$, and
$\Psi_{f,x^0}(y)<0$ in $\mathbb D^n$ if and only if the Lelong number
of $f$ at $x^0$ is strictly positive, otherwise $\Psi_{f,x^0} \equiv
0$.

\begin{defn}
The \emph{local indicator}  $ \Psi_{f,0}$ of a function  $f \in PSH_-(\mathbb
D^n)$ at $x^0=0$ is defined for $y\in \mathbb D^n$ by
$$
\Psi_{f,0}(y)=\underset {R \to +\infty}
{\lim}R^{-1}f[exp(u_k+i\theta_k+R\log|y_k|)].
$$
\end{defn}

This limit exists almost everywhere for $x_k = u_k+i\theta_k$, and does not 
depend on it \cite{lelongandrashkovskii}.

%{\it Remark.} The above limit exists everywhere for
%$x_k=exp(u_k+i\theta_k)$. At  $x^0\not= 0,$ the function
%$\Psi_{f,x^0}(y)$ is defined by means of $f[x^0_k +
%exp(u_k+i\theta_k+R\log|y_k|)].$

Let us fix the system $S:=\{(a_j,\Psi_j)\},1\le j\le N, $ where $a_j\in
\Omega, 1\le j\le N$ and $\Psi_j$ are indicators.

Then by a result of Lelong and Rashkovskii
\cite{lelongandrashkovskii} we have
$$
(dd^c\Psi_j(\cdot - a_j))^n=\tau_j\delta(a_j), 1\le j\le N.
$$

We recall the definition of the generalized Green function due to 
Lelong and
Rashkovskii \cite{lelongandrashkovskii}.

\begin{defn}
$$
G_S(z):= \sup \{u(z): u \in PSH_{-}(\Omega), u(x) \le \Psi_j(x - a_j) + 
C_j,
  \le j \le N\}.
$$
\end{defn}

\begin{rem}
\label{Dirichlet}
If $\Omega$ is a hyperconvex domain in $\mathbb C^n$, then Lelong and
Rashkovskii \cite{lelongandrashkovskii} also showed that the Green
function is the unique solution of the following  Dirichlet problem
(for short we write $G$ instead for $G_S$)

(a) $G \in PSH_- (\Omega) \cap C(\overline{\Omega}); $

(b) $G(z) \to 0$ as $z \to \partial\Omega$;

(c) $\Psi_{G,a_j} = \Psi_j, 1 \le j \le N;$

(d)  $(dd^c G)^n = \sum_{j=1}^N \tau_j \delta (a_j);$
\end{rem}

We now introduce a new generalization of the Lempert function with 
simple
poles, differing from the $\ell_S$ given in the Introduction.

\begin{defn}
\label{newLemp}
\begin{eqnarray*}
\lefteqn{
L_S(z):= \inf \{\sum_{j=1}^N \tau_j\log|\zeta_j|: \exists \varphi \in
\mathcal {O}(\mathbb D, \Omega), \varphi(0)=z, \exists U_j
\mbox{\rm\ a neighborhood of }\zeta_j }
\\
& &
\Psi_j (\varphi(\zeta) -a_j) \le \tau_j \log|\zeta-\zeta_j| + C_j, 
\forall
\zeta \in U_j, 1\le j \le N \}.
\end{eqnarray*}
\end{defn}

Note that for the non-trivial case where $\tau_j \neq 0$, the
conditions imposed on the maps $\varphi$ force $\varphi(\zeta_j)=a_j$.

\begin{lemma} \label{lemma:inequality}
$G_S(z) \le L_S(z)$, for any $z \in \Omega$.
\end{lemma}

\begin{proof}  If $\varphi: \mathbb D \to \Omega $ is an analytic disc
in $\Omega$, with $\varphi(0)=z, \varphi(\zeta_j)=a_j, 1 \le j \le N$
and $\Psi_j \circ \varphi(\zeta) \le \tau_j \log|\zeta-\zeta_j| + C_j,
1\le j \le N,$ then $G\circ \varphi$ is a subharmonic function on
$\mathbb D$, $G \circ \varphi$ is negative and
$$
G_S \circ \varphi(\zeta) \le C_j + \Psi_j \circ \varphi(\zeta)
\le C'_j + \tau_j \log|\zeta-\zeta_j|, \ 1\le j \le N.
$$
Thus $G_S \circ \varphi$
is a member in the defining family for the Green function on $\mathbb
D$ with poles $\zeta_j$ and weights $\tau_j$, and hence,
$$
G_S \circ \varphi(\zeta) \le \sum_{j=1}^N \tau_j \log {\frac {| \zeta_j
- \zeta |}{| 1-\zeta \overline{\zeta_j}|}}.
$$
It implies that
$$
G_S(z)=G_S \circ \varphi(0) \le \sum_{j=1}^N \tau_j \log | \zeta_j| .
$$
Thus $ G_S(z) \le L_S(z), \forall z \in \Omega$.
\end{proof}

We should mention that in general, we don't know how to compare the new
function $L_S$ with the function $\ell_S$ given in the introduction in 
the
case when $\Psi_j (z) = \nu_j \log |z|$ (and therefore $\tau_j = 
\nu_j^n$).

Recall (see e.g. \cite{Garnet}, \cite{Rudin})
that the involutive M\"obius map of $\mathbb D$ which exchanges
  $\xi \in \mathbb D$ and $0$ is given by the following formula:
\begin{equation}
\label{Mobmap}
\phi_{\xi}(\zeta):=\displaystyle{\frac {\xi-\zeta}
{1-\overline{\xi}\zeta}}.
\end{equation}

\begin{lemma}\label{L_S=G_S}

Let $\Omega$ be the polydisc $\mathbb D^n$ in $\mathbb C^n$. If $S$ has
only one pole, and the indicator $\Psi$ is of the following simple kind
$$
\Psi(z)=\max_{1\le j\le n}c_j\log|z_j|,
$$
where the numbers $c_j$ are positive integers, then $L_S(z)=G_S(z),
\forall z\in \mathbb D^n.$
\end{lemma}

\begin{proof} By composing with M\"obius maps in each coordinate, we 
may
reduce ourselves to the case where the pole
$a$ is the origin $0$. By verifying the  Dirichlet problem given by 
Lelong
and Rashkovskii \cite{lelongandrashkovskii}, we have
$$
G_S(z) = \max_{1\le j \le n}c_j\log|z_j|.
$$
We may assume that $\max_{1 \le j \le n}c_j
\log|z_j|=c_{j_0}\log|z_{j_0}|$ for some $1 \le j_0 \le n$.  With this
assumption we have $G_S(z)=c_{j_0}\log|z_{j_0}|$. To prove the Lemma,
  it suffices to show that there exists a mapping $\varphi
\in \mathcal {O}(\mathbb D,\mathbb D^n)$ and $ \zeta_0 \in \mathbb D$
such that

(1) $\varphi (0)=z,$

(2) $\varphi (\zeta_0) = 0, $

(3) $\Psi \circ \varphi (\zeta) \le m \log|\zeta - \zeta_0| +C, \forall
\zeta \in \mathbb D^n,$  where $m:=\prod^n_{j=1}c_j=$ total mass of 
$(dd^c\Psi)^n$

(4) $m\log|\zeta_0| =c_{j_0}\log|z_{j_0}|,$

The condition (3) can be rewritten as follows

(3') $\varphi^{(k)}_j(\zeta_0)=0, 1 \le k \le m_j -1, 1 \le j \le n, $
where $m_j:= m/c_j.$

We fulfill condition (4) by picking $\zeta_0 \in \mathbb D$ such that
$$ |\zeta_0|^{m_{j_0}}= |z_{j_0}|
%^{\displaystyle {\frac {1}{m_{j_0}}}},
$$
and put
$$
\varphi_j(\zeta):=\bigg [\phi_{\zeta_0}(\zeta)\bigg]^{m_j}h_j\bigg
(\phi_{\zeta_0}(\zeta)\bigg), \forall \zeta \in \mathbb D, 1\le j \le
n
$$
where $h_j:\mathbb D \to \overline {\mathbb D}$ is such that
$h_j(\zeta_0)=\displaystyle {\frac {z_j}{{\zeta_0^{m_j}}}}, 1 \le j \le
n.$

Then the function $\varphi =(\varphi_1,\cdots,\varphi_n)$ and $\zeta_0$
satisfy all properties (1), (2), (3') and (4).
\end{proof}

\section{Existence of extremal discs}
\label{extdiscs}

We now extend to this new Lempert function some known properties of its
usual counterpart. The following generalizes \cite[Theorem
2.4, p. 1054]{WikstromAMS}, or in the case of the unit ball 
\cite[Proposition 3,
p. 338]{ThomCCP} (see also \cite[Papers V and VI]{Wikstrom}).

\begin{prop} \label{pro:minusonepoint}
Let $\Omega$ be a convex domain and $S:= \{(a_j,\Psi_j): 1 \le j \le
N\}$ and $S':= \{(a_j,\Psi_j): 1 \le j \le N-1\}$ where $a_j \in
\Omega$ and $\Psi_j$ are indicators centered at $a_j$. Then
$$
L_S (z) \le L_{S'} (z), \mbox{ for all } z \in \Omega.
$$
\end{prop}

The proof of this proposition will be given below in Section 
\ref{addone}.
We will use the shorthand $S' \subset S$ to mean that the sets of poles
are included as noted, and that the indicators remain the same for all
points of the smaller set, as in the above Proposition.

\begin{prop}
\label{convsubset}
Let $\Omega$ be a bounded taut domain, and
$S=\{(a_j,\Psi_j)\}_{j=1,..,N}, N \ge 2.$ If $L_S(z)$ is not
attained by any analytic disc, then
$$L_S(z)\ge
%\underset {S'\varsubsetneq S}
%{\text{min}}
\min_{S'\varsubsetneq S} L_{S'}(z).
$$
In particular, if $\Omega$ is convex and bounded, the conclusion 
becomes
$$
L_S(z)=\min_{S'\varsubsetneq S} L_{S'}(z).
$$
\end{prop}

\begin{proof} The proof of this Proposition is adapted from that of
\cite[Theorem 2.2, p. 1053]{WikstromAMS}.

Take a sequence of analytic discs $\varphi^k$, where
$$
\varphi^k(0)=z \mbox{ and }\Psi_j\circ
\varphi^k(\zeta) \le \tau_j \log|\zeta-\zeta_j^k|+C_j^k , \,
\forall \zeta \in \mathbb D, \, k \ge 1, \, 1 \le j \le
N
$$
such that $\sum_{j=1}^N \tau_j \log|\zeta_j^k|$ converges to
$L_S(z),$ as $k$ tends to $0$.

By passing to a subsequence, using that $\Omega$ is taut, we may
assume that $\varphi^k$ converges locally uniformly to some
$\varphi \in \mathcal{O}(\mathbb D,\Omega).$ Also (if necessary,
by passing to a subsequence again), we may assume that $\zeta_j^k
\to \zeta_j \in \overline{\mathbb D}$, for each $1\le j \le N$, as
$k \to \infty$.

We need to see that for each $\zeta_j\in \mathbb D$,
\begin{equation}
\label{singaj}
\Psi_j\circ
\varphi(\zeta) \le \tau_j \log|\zeta-\zeta_j|+C_j , \,
\mbox{ for } \zeta{ in a neighborhood of }\zeta_j .
\end{equation}

Recall (from \cite{lelongandrashkovskii}) that $\Psi$
being an indicator (centered at $0$) means that
$$ \Psi (z_1, \dots, z_n) = g (\log |z_1|, \dots, \log|z_n|), $$
where $g$ is a convex continuous nonpositive valued function defined on
$(\mathbb R_-)^n$, increasing with respect to each single variable, and
positively homogeneous of degree $1$: $g(\lambda x_1, \dots, \lambda
x_n) = \lambda g( x_1, \dots, x_n)$, for any $\lambda >0$.

We study the situation for a fixed pole $a_j$.
We must have for each $k \ge 0$,
$$
\varphi^k (\zeta^k_j + h) =
(\varphi^k_l (\zeta^k_j + h), 1 \le l \le n)
= ( \alpha_{k,l} h^{m_{k,l}} + O(|h|^{m_{k,l}+1}),1 \le l \le n).
$$
From the above expression,
$$
\Psi_j (\varphi^k (\zeta^k_j + h)) =
g \left( - m_{k,l} + \frac{\log|\alpha_{k,l}|+O(h)}{|\log|h||} \right)
\log |h|,
$$
so the conditions on $\varphi^k$ imply that
\begin{equation}
\label{gless}
g(-m^k) \le \tau_j, \mbox{ where } m^k
:= (m_{k,1}, \dots, m_{k,n}) .
\end{equation}
Passing to a subsequence if needed, we may assume
that $m^k \to m =: (m_1, \dots,m_n) \in (\mathbb N\cup \{\infty\})^n$.
  The uniform convergence on compacta of the
sequence $\varphi^k$ implies that of all derivatives, and that in the
limit $\varphi^{(q)}_l (\zeta_j) = 0$ for $q \le m_l -1$. This, 
together
with (\ref{gless}), proves (\ref{singaj}).

If no $\zeta_j\in \partial \mathbb D$, $\varphi$ is an analytic disc
attaining the infimum in the definition of $L_S(z)$. That is
excluded by our hypothesis. Otherwise, assume after renumbering
the coordinates that $\zeta_j \in 
\mathbb D, 1\le j
\le M$ and $\zeta_j \in \partial \mathbb D$ for $M+1 \le j \le N.$
(Note that not every $\zeta_j$ can be in $\partial\mathbb D,$ as
this would imply that $L_S(z)=0.$) Then $\varphi$ is a member in
the defining family for $L_{S'},$ where
$S':=\{(a_j,\Psi_j)\}_{j=1,..,M},$ and thus $L_S(z)\ge L_{S'}(z).$
\end{proof}

\begin{cor} \label{cor:extremaldisc}
Let $\Omega$ be a bounded taut domain in $\mathbb C^n,$ and let $S$ be
as above. Then for every $z\in \Omega$ there exists an analytic
disc $\varphi,$ such that $\varphi(0)=z,$ passing through a (non
empty) $S_0 \subset S$ such that $\varphi$ attains the infimum in the 
definition
of $L_{S_0}(z)$, and $L_{S_0}(z)=\min_{\emptyset \neq S'
\subset S}L_S(z)$.
\end{cor}

\begin{proof}
If $S$ is a singleton, a normal family argument close to the one used
in the previous proof will show that the corollary is true for this
case.

Otherwise, by the previous proposition, either there is an analytic
attaining the infimum, or $L^S(z)=L_{S_0}(z)$ for some proper subset
$S_0\subset S,$ and $L_{S_0}(z)$ is attained by an analytic disc passing
though $z$ and the points in $S_0$ (otherwise one could pass to a still 
smaller subset).
\end{proof}

As the consequence of Corollary \ref{cor:extremaldisc} and Proposition
\ref{pro:minusonepoint} we have the following.

\begin{thm} \label{thm:infimumattained}
Let $\Omega$ be a bounded convex domain, then
the infimum in the definition of the function $L_S$ is attained by an
extremal disc that passes through a
(non-empty) subset $S' \subset S$ (possibly the whole system $S$).
\end{thm}

However, it would be natural to consider as well the more general case 
of
the relationship between the Lempert functions of two systems
$S:= \{(a_j,\Psi_j): 1 \le j \le
N\}$ and $S':= \{(a_j,\Psi'_j): 1 \le j \le N\}$, where $\Psi_j \le
\Psi'_j$, for any $1 \le j \le N$ ($S'\subset S$ corresponds to the 
case
where the $\Psi'_j$ have $\tau_j= 0$ for $a_j$ outside the pole set of
$S'$). Unfortunately, our generalized Lempert function is not in 
general
monotone when we compare two such generalized pole sets, see a
counter-example below (Proposition \ref{cexnewLemp}). We therefore
introduce a corrected Lempert function $\tilde L$.

\begin{defn} \label{modifLemp}
Let $S:= \{(a_j,\Psi_j): 1 \le j \le N\}$ and $S_1:= \{(a_j,\Psi^1_j):
1 \le j \le N\}$ where $a_j \in \Omega$ and $\Psi_j$, $\Psi^1_j$ are
indicators. %centered at $a_j$.
We define
$$
\tilde L_S (z) := \inf \{ L_{S^1}(z) : \Psi^1_j \ge \Psi_j + C_j, 1
\le j \le N\} .
$$
\end{defn}

\begin{lemma} \label{inegmodifLemp}
$G_S(z) \le \tilde L_S (z) \le L_S (z)$.
\end{lemma}

\begin{proof} The fact that $\tilde L_S (z) \le L_S (z)$ follows from
the definition. For any $S_1$ as in the definition, $L_{S^1}(z)\ge
G_{S^1} (z) \ge G_S(z)$, as follows from Lemma \ref{lemma:inequality}
and the definition of the pluricomplex Green function.
\end{proof}

\section{Adding One Point}
\label{addone}

\begin{proof*}{\it Proof of Proposition \ref{pro:minusonepoint}}

This proof adapts the ideas of
\cite{WikstromAMS} (see also \cite[Proposition 3]{ThomCCP}, 
\cite[Theorem
2.7]{Trao}). Given any $\delta >0$, there exists a holomorphic map
$\varphi$ from the disk to $\Omega$ and points $\zeta_j^0 \in \mathbb
D$, $1\le j \le N-1$, such that
$$
L_{S'} (z) \leq \sum_{j=1}^{N-1} \tau_j \log |\zeta_j^0| \leq L_{S'}
(z) +\delta,
$$
and $\Psi_j \circ \varphi(\zeta) \le \tau_j \log|\zeta-\zeta_j^0| +
C_j$, $1\le j \le N-1$. Let $r<1$ to be specified later. We set
$\varphi^r (\zeta) := \varphi (r\zeta)$. If $r > \max |\zeta_j^0|$,
$1\le j \le N-1$, we then have
$$ \varphi^r (\frac{\zeta_j^0}r) = a_j, 1 \le j \le N-1, $$
and more generally
$$ \Psi_j \circ \varphi^r (\zeta) \le \tau_j \log|r
(\zeta-\frac{\zeta_j^0}r)| + C_j \le \tau_j
\log|(\zeta-\frac{\zeta_j^0}r)| + C_j, 1 \le j \le N-1. $$
We will introduce a correcting term to ensure that the same property
hold for $j=N$, without destroying it for $j \le N-1$. It is no loss of
generality to work in a neighborhood of $0$.

\begin{lemma} \label{locperturb}
Suppose $\Psi$ is an indicator (centered at $0 \in \Omega$) and
$\varphi$ a map from $\mathbb D$ to $\Omega$ such that $\Psi
(\varphi(\zeta)) \le \tau \log|\zeta| + C$. Then there exists
$m=m(\varphi)$ such that for any map $h$ from $\mathbb D$ to $\Omega$
such that $h_j^{(k)} (0)=0$, $1\le j\le n$, $0 \le k \le m$, then $\Psi
( \varphi(\zeta) +h(\zeta)) \le \tau \log|\zeta| + C'$. \end{lemma}

\begin{proof}

The
hypothesis on $\varphi$ implies that $\varphi(0)=0$ and therefore there
exist integers $m_j$ and non-zero complex numbers $\alpha_j$, $1\le j
\le n$, such that $\varphi_j (\zeta) = \alpha_j \zeta^{m_j} (1+o(1))$.
Let $m:= \max_{1\le j \le n} m_j$.
Use the same notations as in the proof of
Proposition \ref{convsubset}.
We then have $$ | \varphi_j (\zeta)
+ h_j (\zeta) | \le M_j | \varphi_j (\zeta) |, $$
for $\zeta$ close enough to $0$. It follows that
\begin{eqnarray*} 
\Psi (\varphi (\zeta) + h (\zeta))& = & g( \log
|\varphi_1 (\zeta) + h_1 (\zeta) |, \dots, \log |\varphi_n (\zeta) +
h_n (\zeta) | )
\\ 
&\le& g( \log |\varphi_1 (\zeta) |+ \log M_1, \dots,
\log |\varphi_n (\zeta)| + \log M_n ) 
\\ 
&=& \left| \log|\zeta| \right|
g( \frac{\log |\varphi_1 (\zeta) |}{\left| \log|\zeta| \right|} +
\frac{\log M_1}{\left| \log|\zeta| \right|}, \dots, \frac{\log
|\varphi_n (\zeta) |}{\left| \log|\zeta| \right|} + \frac{\log
M_n}{\left| \log|\zeta| \right|} ) .
\end{eqnarray*}
Observe that
$$ \lim_{\zeta\to0} \frac{\log |\varphi_j (\zeta) |}{\left| \log|\zeta|
\right|} = - m_j, $$
so that for $|\zeta|$ small enough both the argument of $g$ in the
above formula and $\frac{\log |\varphi_j (\zeta) |}{\left| \log|\zeta|
\right|}$ itself are in a fixed neighborhood of $(-m_1,\dots,-m_n)$,
which we may assume compact within $(\mathbb R_-^*)^n$. From the fact
that the function $g$ is convex, we deduce that it is Lipschitz on any
compact subset of the interior of its domain, therefore for $|\zeta|$
small enough,

\begin{eqnarray*} g( \frac{\log |\varphi_1 (\zeta) |}{\left|
\log|\zeta| \right|} + \frac{\log M_1}{\left| \log|\zeta| \right|},
\dots, \frac{\log |\varphi_n (\zeta) |}{\left| \log|\zeta| \right|} +
\frac{\log M_n}{\left| \log|\zeta| \right|} ) \\ \qquad \le g(
\frac{\log |\varphi_1 (\zeta) |}{\left| \log|\zeta| \right|}, \dots,
\frac{\log |\varphi_n (\zeta) |}{\left| \log|\zeta| \right|} ) + K
\frac{M}{\left| \log|\zeta| \right|},
\end{eqnarray*}
where $K$ is the Lipschitz constant and $M:=\max_j M_j$, therefore

\begin{eqnarray*} g( \log |\varphi_1 (\zeta) + h_1 (\zeta) |, \dots,
\log |\varphi_n (\zeta) + h_n (\zeta) | ) \\ \qquad \qquad \le g( \log
|\varphi_1 (\zeta) |, \dots, \log |\varphi_n (\zeta) | ) + KM = \Psi
(\varphi (\zeta)) + KM,
\end{eqnarray*}
which concludes the proof by the hypothesis on $\varphi$.
\end{proof}

Let $K$ denote the convex hull of $\varphi^r(\overline {\mathbb D})
\cup\{a_N\}$. Since $\varphi^r(\overline {\mathbb D}) \cup \{(a,0)\}
\subset \subset \Omega$, we can find an $\varepsilon>0$ such that the
distance between $K$ and $ \partial \Omega$ is at least $\varepsilon
M_1$ where $M_1 := \sup_{r \overline {\mathbb D}} |a_N- \phi|$.

\begin{lemma} \label{fcnh} Given any $m \in \mathbb N^*$,
there exists $h$ a holomorphic function on
$\mathbb D$ and some $\zeta^* \in \mathbb D$ satisfying
\begin{itemize} \item $h(\mathbb D) \subset U_\varepsilon := \cup_{x
\in [0,1]} D(x,\varepsilon) $, \item $h(0)=0$, \item
$h^{(k)}(\frac{\zeta_j^0}r) =0$, $0 \le k \le m-1$, $1 \le j \le N-1$,
\item $h(\zeta^*) = 1$, and $h^{(k)} (\zeta^*)=0$, $1 \le k \le m-1$.
\end{itemize}
\end{lemma}

Accepting this lemma temporarily, define
$\tilde\varphi(\zeta)=\varphi^r(\zeta)+h(\zeta)(a_N-\varphi^r(\zeta))$.
The definition of $\varepsilon$ and the first condition above show that
$\tilde \varphi(\mathbb D) \subset \Omega$. Clearly,
$\tilde\varphi(0)=z$. Choosing $m$ greater or equal to the maximum of
all the $m$'s that appear in Lemma \ref{locperturb} for various points
$\zeta_j^0/r$, we see that $h(\zeta)(a_N-\varphi^r(\zeta)) =
O((\zeta-\frac{\zeta_j^0}r)^m)$, so that Lemma \ref{locperturb} implies
that for $1 \le j \le N-1$,
$$ \Psi_j \circ \tilde\varphi (\zeta) \le \tau_j
\log|(\zeta-\frac{\zeta_j^0}r)| + C'_j . $$
Finally, one also checks that $\tilde \varphi (\zeta) = a_N +
(h(\zeta)-1)(a_N-\varphi^r(\zeta)) = a_N + O((\zeta-\zeta^*)^m)$, which
for $m$ large enough (depending on $\Psi_N$) will imply \newline
$\Psi_N \circ \tilde\varphi (\zeta) \le \tau_N \log|(\zeta-\zeta^*| +
C_N$. For the mapping $\tilde \varphi$, the logarithmic sum of the
preimages yields
$$ 
\sum_{j=1}^{N-1} \log \left|\frac{\zeta_j^0}r \right| + \log
|\zeta^*| \le \sum_{j=1}^{N-1} \log |\zeta_j^0| + (N-1) \log \frac1r
\le L_{S'} (z) + \delta + (N-1) \log \frac1r. 
$$
Since this construction can be carried out for any $r$ arbitrarily
close to $1$, we have $L_S (z) \le L_{S'} (z)$.
\end{proof*}

\begin{proof*}{{\it Proof of Lemma \ref{fcnh}}} Let $\rho$ be a Riemann
map from $\mathbb D$ to $U_\varepsilon$ so that $\rho(0)=0$. We look
for $h$ under the form $h= \rho \circ h_1$, where $h_1$ is a
holomorphic map from $\mathbb D$ to itself such that
\begin{itemize} \item $h_1(0)=0$, \item $h_1^{(k)}(\frac{\zeta_j^0}r)
=0$, $0 \le k \le m-1$, $1 \le j \le N-1$, and \item there exists
$\zeta^* \in \mathbb D$ such that $h_1 (\zeta^*)= \rho^{-1}(1)$ and
$h_1^{(k)} (\zeta^*)=0$, $1 \le k \le m-1$.
\end{itemize}
The existence of such a function follows from the following 
one-variable Lemma.
%\ref{multPN} in the
%next section.
\end{proof*}

%\section{ A one-variable Lemma}

\begin{lemma}
\label{multPN}
Let $B_0$ be a finite Blaschke product, $\gamma \in \mathbb D$, $m$ a
positive integer. Then there exist $\zeta^* \in \mathbb D$ and $f$
holomorphic from the disk to itself such that
\begin{itemize} \item $f=B_0 \tilde g $, with $\tilde g $ holomorphic
on the disk, \item $f(\zeta^* + h) = \gamma + O(h^m)$.
\end{itemize}
\end{lemma}

\begin{proof} We need to find $g$ holomorphic from the disk to itself
so that the conditions of the lemma are satisfied. We write $\tilde g=g
\circ \phi_{\zeta^*}$, with $\phi_{\zeta^*}$ as in (\ref{Mobmap}). 
Denote
$$ g (z) = \sum_{n \ge 0} a_n z^n, \quad \tilde g (\zeta^* + h)
=\sum_{n \ge 0} \tilde a_n h^n . $$
We have
$$ \phi_{\zeta^*} (\zeta^* + h ) = \frac{-h}{1-|\zeta^*|^2} \frac{1}{1-
\frac{\bar \zeta^* h}{1-|\zeta^*|^2}}, $$
therefore
$$ \phi_{\zeta^*} (\zeta^* + h )^n = \left( \frac{-h}{1-|\zeta^*|^2}
\right)^n \sum_{k \ge 0} \frac1{(n-1)!} (k+1) \dots (k+n-1) \left(
\frac{\bar \zeta^*h}{1-|\zeta^*|^2} \right)^k . $$
By substituting this expression into that of $g$, we find the
relationship between the Taylor coefficients of $g$ and those of
$\tilde g$:
$$ \tilde a_0 = a_0, \quad \tilde a_n = \left( \frac{\bar
\zeta^*}{1-|\zeta^*|^2} \right)^n \sum_{j=1}^n \frac{(n-j+1)(n-j+2)
\dots (n-1)}{(j-1)!} (\bar \zeta^*)^{-j}a_j . $$
Solving this triangular system of linear equations, we see that there
exist coefficients $c(j,n)$, independent of $\zeta^*$, such that
$$ a_j = (\bar \zeta^*)^{j} \sum_{n=1}^j \left(
\frac{1-|\zeta^*|^2}{\bar \zeta^*} \right)^n c(j,n) \tilde a_n . $$
Now let us write the necessary and sufficient conditions on $\tilde
a_n$ to get $f(\zeta^* + h) = \gamma + O(h^m)$. We have
$$ B_0 (\zeta^* + h) \tilde g (\zeta^* + h) = \sum_{k \ge 0} \left(
\sum_{j=0}^k \frac{B_0^{(j)} (\zeta^*)}{j!} \tilde a_{k-j} \right) h^k
, $$
so in order to get the required local expansion, we must have $\tilde
a_0 = \gamma / B_0 (\zeta^*)$ and recursively
$$ \tilde a_k = -\frac1{B_0 (\zeta^*)} \sum_{j=1}^k \frac{B_0^{(j)}
(\zeta^*)}{j!} \tilde a_{k-j} , \quad 1 \le k \le m-1. $$
By induction, we see that there exist polynomials in $k+1$ variables
with coefficients depending only on $k$, $P_k (X_0, \dots, X_k)$ such
that those relations are equivalent to
$$ \tilde a_k = \frac1{B_0 (\zeta^*)^{k+1}} P_k (B_0 (\zeta^*), B'_0
(\zeta^*), \dots, B^{(k)}_0 (\zeta^*)). $$
Finally, for a given $\zeta^*$, a function $\tilde g$ satisfying the
requirements of the Lemma will exist if and only if we can find a
function $g$ holomorphic from the unit disk to itself such that $a_0 =
g(\zeta^*)$, and for $1 \le j \le m-1$,
$$ a_j = (\bar \zeta^*)^{j} \sum_{n=1}^j \left(
\frac{1-|\zeta^*|^2}{\bar \zeta^*} \right)^n c(j,n) \frac1{B_0
(\zeta^*)^{n+1}} P_n (B_0 (\zeta^*), B'_0 (\zeta^*), \dots, B^{(n)}_0
(\zeta^*)). $$
Since $B_0$ is a finite Blaschke product, we know that there exist some
$\eta >0$ so that it is holomorphic in a neighborhood of the compact
annulus $\{ 1-\eta \le |z| \le 1+\eta \}$, and that all $|B_0^{(j)}|$
are bounded above, and $|B_0|$ is bounded and bounded away from zero on
this annulus. Choose $\zeta^*$ within this annulus. Since $\gamma \in
\mathbb D$, we may choose $|\zeta^*|$ close enough to $1$ so that
$|\tilde a_0 | < 1$, and we see that $|a_j| \le C_j (1-|\zeta^*|^2)$,
for $1 \le j \le m-1$. Since we can choose $(1-|\zeta^*|^2)$ as small
as we need, the proof will conclude with the next Lemma.
\end{proof}

\begin{lemma} Given any $r<1$ and an integer $m \ge 1$, there exists
$\varepsilon = \varepsilon (r,m)$ so that for any $\beta$ such that
$|\beta|\le r$ and $a_j$, $1\le j \le m-1$ such that $|a_j| \le
\varepsilon$, $1\le j \le m-1$, then there exists a holomorphic
function $g$ from the disk to itself such that
$$ g(z) = \beta + a_1 z + \cdots + a_{m-1} z^{m-1} + O(z^m) . $$
\end{lemma}

\begin{proof} First it is clear that we may reduce ourselves to $\beta
=r$, for if we assume the problem solved for $\beta =r$, we can pick $g
= (\beta/r) g_1$, with
$$ g_1(z) = r + \frac{r}{\beta} \left( a_1 z + \cdots + a_{m-1} z^{m-1}
\right) + O(z^m) . $$
We proceed by induction, in the spirit of the proof of Carath\'eodory's
theorem about approximation by Blaschke products \cite[p. 6, Theorem
2.1]{Garnet}. For $m=0$ the constant function equal to $r$ will do.
Suppose the property is known for $m-1$. For any function $g$
holomorphic on the disk with $g(0)=r$, define a new function $g_2$ by
$$ g_2 (z) := \frac1z \frac{g(z)-r}{1- r g(z)}, \mbox{ i.e. } g(z) =
\frac{zg_2(z) + r}{1+ r z g_2(z)}. $$
Then $g$ sends the disk to the disk if and only if $g_2$ does, and the
Taylor coefficients of $g$ up to order $m-1$ are determined by the
Taylor coefficients of $g_2$ up to order $m-2$ and vice-versa. Indeed,
suppose that $g_3 = g_2 + f$, with $f(z)=O(z^{m-1})$, then
$$ \frac{zg_3(z) + r}{1+ r z g_3(z)} - \frac{zg_2(z) + r}{1+ r z
g_2(z)} =\frac{zf(z)(zg_2(z) +r - 1 - rz g_2(z))}{(1+ r z g_3(z))(1+ r
z g_2(z))} = O(z^{m}). $$
Conversely, if we are given $a_j$, $1\le j \le m-1$ and $g(z) = \beta +
a_1 z + \cdots + a_{m-1} z^{m-1} + O(z^m) $,
$$ \frac1z \frac{g(z)-r}{1- r g(z)} = \frac{g(z)-r}{1-r^2} \sum_{k\ge
0} \frac{r^k (g(z)-r)^k}{(1-r^2)^k} ; $$
by expanding out and collecting terms, we see that there exist
polynomials \newline $q_{j,r} (a_1,\dots,a_{j})$, $1 \le j \le m-1$,
such that $q_{1,r} (a_1) = (1-r^2)^{-1} a_1$, $q_j (0,\dots,0)=0$ for
any $j$, and
$$ zg_2(z) = \sum_{j=1}^{m-1} q_{j,r} (a_1,\dots,a_{j}) z^j + O(z^m),
$$ 
therefore
$$ g_2(z) = (1-r^2)^{-1} a_1 + \sum_{j=1}^{m-2} q_{j+1,r}
(a_1,\dots,a_{j}) z^j + O(z^{m-1}) . $$
We pick $\varepsilon (r, m) \le r (1-r^2)$, so that the first term 
$q_{1,r}(a_1)$ is
less than $r$ in modulus ; then by continuity of each $q_j$ one can
choose $\varepsilon (r, m)$ small enough so that when $|a_j| \le
\varepsilon (r, m)$, $2\le j \le m-1$, we have
$$ |q_{j,r} (a_1,\dots,a_{j})| \le \varepsilon (r, m-1), $$
so that the result for $m-1$ yields the existence of $g_2$, therefore
that of $g_1$.
\end{proof}

\section{ Examples in the bidisc}
\label{bidisc}

First, we would like to give one case where the Green function with 
several poles
and indicator singularities is equal to its generalized Lempert
counterpart. This is analogous in spirit to the result of Carlehed and
Wiegerinck about
the Green function with several poles in the bidisc \cite{pointsext},
\cite{CarlWieg} (but easier).

\begin{thm}
\label{thm:G=LinthecaseNpoles}
Let $ \Psi_m (z) = \max\, \{m \log |z_1 |;\log|z_2| \}$, for any $m \in
\mathbb N^*$.
\newline
Let $a_1, a_2,\dots , a_N \in \mathbb D$, and
$$
S := \{ ((a_1,0); \Psi_{m_1}); \dots, ((a_N,0); \Psi_{m_N}) \}.
$$
  Then for any $z\in \mathbb D^2$,
$$
L_S(z)=G_S(z)= \max \{ \sum_{j=1}^N m_j \log \phi_{a_j} (z_1) ; \log
|z_2| \}.
$$

As a consequence, if $a_{j,i}^{(k)} \in \mathbb D$, $1\le j \le N$, 
$1\le i
\le m_j$, are distinct points which verify
$$
\lim_{k\to\infty} a_{j,i}^{(k)} = a_j, \quad, 1\le i \le m_j,
$$
and $S^{(k)}$ the pole system made up of all the  $a_{j,i}^{(k)}$ with
equal weight $1$,
\newline
then $\lim_{k\to\infty} L_{S^{(k)}} (z) = L_S(z)$
  and $\lim_{k\to\infty} G_{S^{(k)}} (z) = G_S(z)$, for
any $z \in \mathbb D^2$.
\end{thm}

\begin{proof} %First of all, we would like to remark that, with the
%hypothesis, the new Green (resp. Lempert) function is the same as the
%usual Green (resp. Lempert) function with $N$ poles $(a_j,0), 1\le j
%\le N$ and weights are equal.

First of all, the Green function has the formula given above.
To prove this assertion it suffices to show that the function defined
by the right hand side verifies the Dirichlet problem
in Remark \ref{Dirichlet}. Indeed the conditions (a), (b) and (c)
are trivially fulfilled. The last condition follows from the following
theorem of Zeriahi \cite{Zeriahi}.

\noindent {\bf Theorem.} {\it For $i=1,2,$ let $\Omega_i$ be an open
set in $\mathbb C^{n_i},$ and $u_i$ a locally bounded plurisubharmonic
function in $\Omega_i,$ such that $(dd^c u_i)^{n_i}=0$ in $\Omega_i.$
Define $v(z_1,z_2)=\max\{u_1(z_1),u_2(z_2)\}, n=n_1+n_2.$ Then
$(dd^c v)^n=0.$ }

By our definition,
\begin{multline*}
L_S(z)= \inf\{ \sum_{j=1}^N m_j \log|\zeta_j| : \exists \varphi \in
\mathcal {O}(\mathbb D, \mathbb D^2),
\\
\varphi(0)=z,
\varphi_1(\zeta_j)=a_j, \varphi_2^{(k)} (\zeta_j)=0,\, 0 \le k \le 
m_j-1, 1\le j \le N \} .
\end{multline*}

If $z_1 \in \{a_1, \dots, a_n\}$, say $z_1 =a_1$, then picking
$\zeta_1^{m_1} = z_2$ and $\varphi(\zeta) = (a_1, \zeta^{m_1})$, we see 
by
Proposition \ref{pro:minusonepoint}
that
$$
\log |z_2| = m_1 \log |\zeta_2| \ge L_{((a_1,0),\Psi_{m_1})} (z) \ge
L_S(z) \ge G_S(z) = \log |z_2|,
$$
so there is equality throughout.

If $z_1 \notin \{a_1, \dots, a_n\}$, we may reduce ourselves
to $z = (0,\gamma)$ and $|a_1| \ge |a_2| \ge ...
\ge |a_N| >0$. Then
$$
G_S(z)= \max \{ \log |a_1^{m_1} \cdot a_2^{m_2}\cdots a_N^{m_N}|; 
\log|\gamma|\}.
$$

  We will use induction on $N$. When $N=1$ the equality follows from
Lemma \ref{L_S=G_S}. Suppose that $N>1$ and the theorem is proved for
$N-1$. We consider three cases.

{\it Case 1.} $|\gamma| \le |a_1^{m_1} \cdot a_2^{m_2}\cdots 
a_N^{m_N}|$.

Then $G_S(z) = \log |a_1^{m_1} \cdot a_2^{m_2}\cdots a_N^{m_N}|$. The map
$$
\zeta \mapsto
\biggl (\zeta,
\displaystyle {\frac {\gamma}
{a_1^{m_1} \cdot a_2^{m_2}\cdots a_N^{m_N}}}
\prod^N_{j=1} \displaystyle {\biggl ({\frac {a_j -
\zeta} {1- \overline {a_j} \zeta}}\biggr )^{m_j}}\biggr )
$$
verifies all the requirements with $\zeta_j =a_j$.
This implies that $G(z)=L(z).$

{\it Case 2.} $|\gamma| \ge |a_2^{m_2}\cdots a_N^{m_N}|$.

Then $G(z) =
\log|\gamma|$. Moreover, $G(z)$ is also equal to the Green function
$G_1(z)$ for the system with $N-1$ poles
$$
S_1 := \{ ((a_2,0); \Psi_{m_2}); \dots, ((a_N,0); \Psi_{m_N}) \}.
$$
By induction, $G_1 = L_1,$ where $L_1$ is the generalized
Lempert function with respect to $S_1$.
On the other hand, we always have $L_S(z)\le
L_1(z)$ by Proposition \ref{pro:minusonepoint}. Hence $G_S(z) = L_S(z).$

{\it Case 3.} $|a_1^{m_1} \cdot a_2^{m_2}\cdots a_N^{m_N}| < |\gamma| 
< |a_2^{m_2}\cdots a_N^{m_N}|.$

We now
show that the $G_S(z)=\log|\gamma|$ is also equal to the new Lempert
function, and the infimum in the definition of the new Lempert function
is attained by an extremal disc $\varphi$  passing through all poles
$(a_1,0); (a_2,0); ...; (a_N,0)$ and $z$.

Set $M:= \sum_{j=1}^N m_j$ and define $r \in (0,1)$ by
$$
r = \root M \of {\displaystyle {\frac
{|a_1^{m_1} \cdot a_2^{m_2}\cdots a_N^{m_N}|} {|\gamma|}}}.
$$
We have, for any $1 \le j \le N$,
$$
|a_j|^M < |a_j|^{m_1} \le |a_1|^{m_1} < r^M
$$
by the hypothesis on $\gamma$. So
$a_j/r \in \mathbb D$. We
introduce the map  $\varphi : \mathbb D \mapsto \mathbb D^2$ given by
$$
\varphi (\zeta) = \biggl (r \zeta ,
e^{i\theta}\prod^N_{j=1}\displaystyle {\biggl ( \frac {\zeta_j - \zeta}
{1-\overline {\zeta_j}\zeta}} \biggr )^{m_j} \biggr ),
$$
  where $\zeta_j = \displaystyle {\frac {a_j} r}$, $1 \le j \le N$,
and $\theta$ is chosen such that
$$
e^{i\theta}  \left({\frac {a_1} r}\right)^{m_1}
  \left({\frac {a_2}  r}\right)^{m_2} \cdots
\left({\frac {a_N} r}\right)^{m_N} = \gamma.
$$
  It is easy to verify that $\varphi$ verifies the conditions in the
  definition of $L_S$ and that
$|\zeta_1^{m_1} \cdot \zeta_2^{m_2} \cdots \zeta_N^{m_N}| = |\gamma|$. Hence, $\varphi$ is an
extremal disc for the new Lempert function, and $G_S(z) = L_S(z)$ in this
case.
\end{proof}

We will now give some negative results,
mainly that the generalized Green function can be different from the
generalized
Lempert function as given in Definition \ref{newLemp}.

We shall need some
notation, to be used in this section and the next one.

For $z \in \mathbb D^2$, we will use the following indicators:
\begin{multline}
\Psi_0 (z) := \max(\log|z_1|,\log|z_2|), \\
\Psi_H (z) := \max(2\log|z_1|,\log|z_2|) ,
\quad \Psi_V (z) := \max(\log|z_1|,2\log|z_2|).
\end{multline}
Here $H$ stands for "horizontal" and $V$ for "vertical", for the 
obvious
reasons : for $a \in \mathbb D^2$, 
\newline
$\Psi_j ( \varphi (\zeta) -a) \le \tau_j \log|\zeta -\zeta_0| + C$
translates to ($\tau_0=1$, $\tau_H=\tau_V=2$):
\begin{eqnarray*}
\varphi (\zeta_0) = a, &\mbox{ when }& j =0, \\
\varphi (\zeta_0) = a, \varphi'_2(\zeta_0) =0 &\mbox{ when }& j =H, \\
\varphi (\zeta_0) = a, \varphi'_1(\zeta_0) =0 &\mbox{ when }& j =V .
\end{eqnarray*}
For $a$, $b \in \mathbb D$, let
\begin{eqnarray*}
S_{a0} &:=& \{ ((a,0), \Psi_0)\} = \{ (a,0) \} \\
S_{a0b0} &:=& \{ ((a,0), \Psi_0);((b,0), \Psi_0) \} = \{ (a,0);(b,0) \} 
\\
S_{aV} &:=& \{ ((a,0), \Psi_V) \} \\
S_{bV} &:=& \{ ((a,0), \Psi_V) \} \\
S_{a0bV} &:=& \{ ((a,0), \Psi_0);((b,0), \Psi_V) \} \\
S_{aVbV} &:=& \{ ((a,0), \Psi_V);((b,0), \Psi_V) \}.
\end{eqnarray*}
We will denote with the corresponding subscripts the pertinent Green 
and
Lempert functions, e.g. $G_{a0bV}$, $L_{a0bV}$, $\tilde L_{a0bV}$, etc.
A special case of Theorem \ref{thm:G=LinthecaseNpoles} is that 
$L_{aHb0} = G_{aHb0}$ for any $a$
and $b$ in the disc, for instance.

We start by giving an example of a situation where $\tilde L_S(z) < 
L_S(z)$,
with $S=S_{a0bV}$.

\begin{prop}
\label{cexnewLemp}
For $z_1\in \mathbb D$, $L_{a0bV}(z_1,0) > L_{a0b0}(z_1,0)$, and 
therefore
\newline
$L_{a0bV}(z_1,0) > \tilde L_{a0bV} (z_1,0) \ge G_{a0bV}(z_1,0)$.
\end{prop}

\begin{proof}
By the above, $L_{a0b0}(z_1,0) = G_{a0b0}(z_1,0) = \log|\phi_a(z_1)|+
\log|\phi_b(z_1)|$, where $\phi_a$ and $\phi_b$ are as in 
(\ref{Mobmap}). We have
\begin{multline*}
L_{a0bV}(z_1,0)= \inf\{\log|\zeta_1|+2\log|\zeta_2|: \exists \varphi 
\in
\mathcal {O}(\mathbb D, \mathbb D^2),
\\
\varphi(0)=(z_1,0),
\varphi(\zeta_1)=(a,0), \varphi (\zeta_2)=(b,0)\; \text {and} \;
\varphi_1'(\zeta_2)=0\} ,
\end{multline*}
\begin{multline*}
L_{a0}(z_1,0)= \inf\{\log|\zeta_1| : \exists \varphi \in
\mathcal {O}(\mathbb D, \mathbb D^2),
\varphi(0)=(z_1,0) \text { and }
\varphi(\zeta_1)=(a,0) \}, \\
L_{bV}(z_1,0)= \inf\{ 2\log|\zeta_2|: \exists \varphi \in
\mathcal {O}(\mathbb D, \mathbb D^2), \hfill \\
\varphi(0)=(z_1,0),
\varphi (\zeta_2)=(b,0) \mbox{ and }
\varphi_1'(\zeta_2)=0\}.
\end{multline*}
So $L_{a0bV}(z_1,0) \ge L_{a0}(z_1,0) + L_{bV}(z_1,0)$, since each of 
the
infima on the right hand side is taken over a  family of maps $\varphi$
which is wider than the one used in the definition of $L_{a0bV}$.

By Lemma \ref{L_S=G_S}, $L_{a0}(z_1,0) = \log |\phi_a(z_1)|$,
$L_{bV}(z_1,0) = \log|\phi_b(z_1)|$.

Now suppose that $L_{a0bV}(z_1,0) \le L_{a0b0}(z_1,0)$. This means
$$
L_{a0bV}(z_1,0) \le G_{a0b0}(z_1,0) =  L_{a0}(z_1,0) + L_{bV}(z_1,0),
$$
so there is equality throughout. Since $ L_{a0bV}(z_1,0) < 
\min(L_{a0}(z_1,0),
L_{bV}(z_1,0))$, Proposition \ref{convsubset} shows that the infimum in 
the
definition of $ L_{a0bV}$ is attained by a map $\varphi$. It follows
from the Schwarz Lemma applied to $a$ and $z_1$ that its
first coordinate $\varphi_1$ is a M\"obius map of the disc.  But we 
also
had to have $\varphi_1'(\zeta_2) =0$. This is a contradiction.
\end{proof}

The following example is similar, and will be useful in the final
construction.

\begin{prop}
\label{prop:G2ab<L2ab}
If $a \neq b \in \mathbb D$ and %$b=-a$ and
$|\gamma|^2<|ab|$, then
$$
G_{aVbV}(0,\gamma)<L_{aVbV}(0,\gamma).
$$
\end{prop}

\begin{proof} First of all we can rewrite the generalized Lempert 
function as
follows
\begin{multline*}
L_{aVbV}(z)= \inf\{2\log|\zeta_1|+2\log|\zeta_2|: \exists \varphi \in
\mathcal {O}(\mathbb D, \mathbb D^2), \varphi(0)=z, \\
\varphi(\zeta_1)=(a,0), \varphi (\zeta_2)=(b,0)\; \text {and} \;
\varphi_1'(\zeta_1)=0,\varphi_1'(\zeta_2)=0\} .
\end{multline*}
As in the proof of Proposition \ref{cexnewLemp},
by Lemma \ref{L_S=G_S} we have
$$
L_{aV}(z)=G_{aV}(z)=\max \{\log|\phi_a(z_1)|;2\log|z_2|\}, \quad 
\forall z \in \mathbb
D^2,
$$
and similarly for $L_{bV}(z)=G_{bV}(z)$.

By using the  Dirichlet problem given by Lelong and
Rashkovskii \cite{lelongandrashkovskii}, we can verify
that
$$
G_{aVbV}(z) = \max \{\log|\phi_a(z_1)| + \log|\phi_b(z_1)|; 2\log|z_2| 
\}.
$$
%Without loss of generality we may assume that $a>0$ (if need, we can
%use a
%rotation).
Since $|\gamma|^2<|ab|$, $G_{aVbV}(0,\gamma)= \log |a| + \log |b|$.

From
Lemma \ref{lemma:inequality} we already know
$G_{aVbV}(z) \le L_{aVbV}(z)$, for any $z\in \mathbb D^2$. Suppose
equality holds at $z_0:=(0,\gamma)$. Then, by using Lemma \ref{L_S=G_S} 
and the
definition
of $L_{aVbV}$ we have
\begin{multline*}
G_{aVbV}(z_0) =
\log|a|+\log|b| \le
\\
G_{aV}(z_0) + G_{bV}(z_0) = L_{aV}(z_0) + L_{bV}(z_0) \le
L_{aVbV}(z_0) = G_{aVbV}(z_0).
\end{multline*}
Hence equality would hold throughout.
Now, by Proposition \ref{convsubset},
the infimum in the definition of $L_{aVbV}$ is
attained by an extremal disc $\varphi$ that passes through both $(a,0)$
and $(b,0)$. It follows that $\varphi$ must be extremal for
$L_{aV}$ and $L_{bV}$. We will prove that this is impossible.

First of all we characterize all extremal discs for $L_{aV}$. Let
$\varphi = (\varphi_1,\varphi_2)$ be such a disc. By the
definition there exists $\zeta_1 \in \mathbb D$ such that
$\varphi(0)=(0,\gamma)$, $\varphi(\zeta_1)=(a,0)$,
$\varphi'_1(\zeta_1)=0$, $|\zeta_1|^2=a$.

Setting $g:=\phi_a\circ\varphi_1\circ \phi_{\zeta_1}$, we have
$$
g(0)=0, g'(0)=0, g(\zeta_1)=a, |\zeta_1|^2=a.
$$
The Schwarz Lemma
now gives $g(\zeta)=e^{i\theta}\zeta^2$, where $\theta \in \mathbb 
R$. It implies that
$$
\varphi_1(\zeta)=\phi_a\bigg (e^{i\theta}\bigg
(\phi_{\zeta_1}(\zeta)\bigg )^2 \bigg ), \quad \forall \zeta \in 
\mathbb
D.
$$
If the function $\varphi$ is an extremal disc for $L_{bV}$, then
there is $\zeta_2 \in \mathbb D$ such that
$$
\varphi_1(0)=0,
\varphi_1(\zeta_2)=b, \varphi'_1(\zeta_2)=0,|\zeta_2|^2=a.
$$
Clearly $\zeta_1 \not= \zeta_2$ since $a \not= b$. Since $\varphi_1$ 
only
has one critical point, the condition $\varphi_1'(\zeta_2)=0$ is
not verified, so we have a contradiction.
\end{proof}

\begin{prop}
\label{prop:G2ab<L2b}
If $a \neq b \in \mathbb D$,  $|\gamma|<|a|$, and
$|\gamma|^2<|ab|$, then
$$
G_{aVbV}(0,\gamma)<L_{a0bV}(0,\gamma).
$$
\end{prop}

\begin{proof}
The arguments are similar to those in the proof of the above 
proposition,
so we only indicate the differences. As in the proof of Proposition
\ref{cexnewLemp},
$L_{a0bV}(z) \ge L_{a0}(z) +  L_{bV}(z) = G_{a0}(z) +  G_{bV}(z)  $
by Lemma \ref{L_S=G_S} ; because of the value of
$|\gamma|$, this is equal to $G_{aVbV(z)}$. So if the conclusion was 
not
true, equality would have to hold throughout, but the extremal disc
$\varphi$ in the definition of $L_{a0}(0,\gamma)$ would have to have a
M\"obius map for its first coordinate $\varphi_1$, and since this has 
no
critical point, it could not be extremal for $L_{bV}(0,\gamma)$.
\end{proof}

We want to see that, in some cases, $\tilde L$ is a better candidate for 
the limit of the usual Lempert functions when poles coalesce.

Let $\varepsilon \in \mathbb C$. Define $L^\varepsilon =
L^\varepsilon_{a0bV}$ to be
the
usual Lempert function with three poles $(a,0); (b,0);
(b,\varepsilon)\in \mathbb D^2 $, and
$\tilde L_{a0bV}$ as above, here
$\tilde L_{a0bV}= \min \{ L_{a0bV}, L_{a0b0} \}$.  %\tilde L_{a0bV}

Denote also $G_{aVbV}^\varepsilon$ (resp. $L_{aVbV}^\varepsilon$) the
usual Green (resp. usual Lempert) function with  four poles
$\{(a,0);(b,0);(b,\varepsilon);(a,\varepsilon)\}$. By using the product
property of the Green fuction \cite{Edigarian},
$$
G_{aVbV}^\varepsilon(z) =\max\bigg \{  \log \bigg |
\displaystyle {\frac {a-z_1}{1-z_1\overline{a}}}\bigg |+\log \bigg |
\displaystyle {\frac {b-z_1}{1-z_1\overline{b}}}\bigg |;\log|z_2|+\log
\bigg | \displaystyle {\frac {\varepsilon - z_2}
{1-z_2\overline{\varepsilon}}}\bigg | \bigg \} .
$$

\begin{thm}
\label{liml3bepsilon=}
$
\lim_{\varepsilon \to
\infty}L^\varepsilon_{a0bV} (z)=\tilde L_{a0bV}(z)$, for any $ z \in
\mathbb D^2$.
\end{thm}

\begin{proof} We consider two cases

{\it Case 1.} $\log|z_2| \le \log \bigg | \displaystyle {\frac {a-z_1}
{1 - z_1\overline{a}}}\bigg | + \log \bigg | \displaystyle {\frac
{b-z_1}{1 - z_1\overline{b}}}\bigg |.$

By the hypothesis and Lemma \ref{L_S=G_S}, we have
$$
L_{a0b0}(z) =\log \bigg | \displaystyle {\frac {a-z_1} {1 -
z_1\overline{a}}}\bigg | + \log \bigg | \displaystyle {\frac {b-z_1} {1
- z_1\overline{b}}}\bigg |
=
L_{a0}(z)+L_{b0}(z) \le L_{a0bV}(z).
$$
By Theorem
\ref{thm:G=LinthecaseNpoles}, Lemma \ref{lemma:inequality}, and the 
usual
inequalities between Lempert functions \cite[Theorem
2.4]{WikstromAMS} or Proposition \ref{pro:minusonepoint} above
for simple poles, we have
\begin{multline*}
\log \bigg | \displaystyle {\frac {a-z_1} {1-z_1\overline{a}}}\bigg|
+  \log \bigg | \displaystyle {\frac {b-z_1} {1-z_1\overline{b}}}\bigg|
=
L_{a0b0}(z) \ge L_{a0bV}^\varepsilon(z)
\ge
L_{aVbV}^\varepsilon(z)
\ge
G_{aVbV}^\varepsilon(z) \\
= \log \bigg | \displaystyle {\frac {a-z_1} {1 -
z_1\overline{a}}}\bigg | + \log \bigg | \displaystyle {\frac {b-z_1} {1
- z_1\overline{b}}}\bigg | .
\end{multline*}
Thus,
$$
L_{a0bV}^\varepsilon(z)=L_{a0b0}(z)
=\log \bigg | \displaystyle {\frac
{a-z_1} {1-z_1\overline{a}}}\bigg |
+
\log \bigg |
\displaystyle {\frac
{b-z_1} {1-z_1\overline{b}}}\bigg |.
$$
Since $L_{a0b0}(z)\le L_{a0bV}(z)$ in this case, we have $\underset
{\varepsilon \to 0}{\lim}L_{a0bV}^\varepsilon=L_{a0b0}=\tilde 
L_{a0bV}$.

{\it Case 2.} $\log|z_2| > \log \bigg | \displaystyle {\frac {a-z_1} {1
- z_1\overline{a}}}\bigg | + \log \bigg | \displaystyle {\frac {b-z_1}
{1 - z_1\overline{b}}}\bigg |.$ Then $L_{a0b0}(z)=\log|z_2|$. We now
divide the proof in two steps. In Lemma \ref{limsup} we prove the 
inequality
$\underset {\varepsilon \to 0}{\limsup}\; L_{a0bV}^\varepsilon \le
\tilde L_{a0bV}$. In Lemma \ref{liminf}  we show that  $\underset
{\varepsilon \to 0}{\liminf} \;L_{a0bV}^\varepsilon \ge \tilde 
L_{a0bV}.$
Then we have $\underset {\varepsilon \to 0}{\text{lim }}\;
L_{a0bV}^\varepsilon =\tilde L_{a0bV}$.

\begin{lemma}
\label{limsup}
For all $z \in \mathbb D^2$ such that $\log|z_2| > \log \bigg | 
\displaystyle {\frac {a-z_1} {1
- z_1\overline{a}}}\bigg | + \log \bigg | \displaystyle {\frac {b-z_1}
{1 - z_1\overline{b}}}\bigg |$, then
$\underset {\varepsilon \to
0}{\limsup}\; L_{a0bV}^\varepsilon \le \tilde L_{a0bV}$.
\end{lemma}

\begin{proof}

  By the monotonicity property of Coman's Lempert function 
\cite[Theorem
2.4]{WikstromAMS} we have $L_{a0bV}^\varepsilon \le L_{a0b0}$ , and
hence $\underset {\varepsilon \to 0}{\limsup}\;L_{a0bV}^\varepsilon \le
L_{a0b0}$. Thus we only need to prove that $\underset {\varepsilon \to
0}{\limsup}L_{a0bV}^\varepsilon\le L_{a0bV}.$

Let $\varphi \in \mathcal {O}(\mathbb D,\mathbb D^2)$ be an  analytic
disc and $\zeta_1, \zeta_2 \in \mathbb D$ such that
\begin{equation}
\label{candidate}
\varphi(0)=z, \;
\varphi(\zeta_1)=(a,0), \; \varphi(\zeta_2)=(b,0), \; 
\varphi'_1(\zeta_2)=0.
\end{equation}
We now consider two cases

$\bullet$ If  $\varphi'_2(\zeta_2)=0$ then $\varphi_2 \in
\mathcal{O}(\mathbb D,\mathbb D)$  satisfies $\varphi_2(0)=z_2,
\varphi_2(\zeta_1)=\varphi_2(\zeta_2)=\varphi_2'(\zeta_2)=0.$ Thus
$$\varphi_2(\zeta)=\bigg( \displaystyle {\frac {\zeta_1-\zeta}{1-\zeta
\overline{\zeta_1}}}\bigg)\bigg( \displaystyle {\frac
{\zeta_2-\zeta}{1-\zeta \overline{\zeta_2}}}\bigg)^2h(\zeta),$$
where $h:\mathbb D \to \overline{\mathbb D}.$ The equality
$\varphi_2(0)=z_2$ implies that $h(0)=\displaystyle{\frac {z_2}
{\zeta_1\zeta_2^2}}\in \overline{\mathbb D}.$ Hence
$L_{a0b0}(z)=\log|z_2| \le \log|\zeta_1|+2\log|\zeta_2|.$

Using the estimate $\underset {\varepsilon \to
0}{\limsup}\;L_{a0bV}^\varepsilon(z) \le L_{a0b0}(z), \forall z \in
\mathbb D^2,$ we get $\underset {\varepsilon \to
0}{\limsup}\;L_{a0bV}^\varepsilon \le \log|\zeta_1|+2\log|\zeta_2|.$

$\bullet$  If $\varphi'_2(\zeta_2)\not= 0.$ Let $\varepsilon \in
\mathbb C$ be such that $|\varepsilon|$ small enough. We will show that
there exist
$\widetilde{\varphi}\in \mathcal{O}(\mathbb D,\mathbb D^2),
\widetilde{\zeta_j} \in \mathbb D, 1 \le j \le 3$ such that
$$
\widetilde{\varphi}(0)=z,
\widetilde{\varphi}(\widetilde{\zeta_1})=(a,0),
\widetilde{\varphi}(\widetilde{\zeta_2})=(b,0),
\widetilde{\varphi}(\widetilde{\zeta_3})=(b,\varepsilon)
$$
and
$$\sum_{j=1}^3\log|\widetilde{\zeta_j}| \to
\log|\zeta_1|+2\log|\zeta_2| \,\text{as}\, \varepsilon \to 0.$$
This will show that $\limsup_{\varepsilon\to 0} L_{a0bV}^\varepsilon(z)\le
\log|\zeta_1|+2\log|\zeta_2|$.

Let $0<r<1$ to be chosen later. Define $\varphi^r:\mathbb D \to \mathbb
D^2$ by $\varphi^r(\zeta)=\varphi(r\zeta), \forall \zeta \in \mathbb
D.$ It is  easy to see that $\varphi^r(0)=z,
\varphi^r(\zeta_1/r)=(a,0), \varphi^r(\zeta_2/r)=(b,0)$. The Invariant
Schwarz Lemma (see \cite[Chap. I, Lemma 1.2, p. 2]{Garnet}) says that 
if
$f$ is a holomorphic function from the unit disc to itself, then
$d_G(f(z),f(w)) \le d_G(z,w)$, where $d_G(z,w):=|\phi_z(w)|$, $\phi_z$
being defined as in (\ref{Mobmap}). Then
$(\varphi_j^r(\zeta),z_j)=d_G(\varphi_j(r\zeta),\varphi_j(0))\le
d_G(r\zeta,0)\le r$, for all $\zeta \in \mathbb D$, $j=1,2$. This implies 
that
$$
|\varphi_j^r(\zeta)|\le \displaystyle{\frac {r+|z_j|}  {1+r|z_j|}}
\le 1-s_j (1-r),\quad \forall \zeta \in \mathbb D, \quad j=1,2,
$$
where $s_j=\bigg ( \displaystyle{\frac {1-|z_j|} {1+|z_j|}}\bigg ),
j=1,2,$ depend only on $z$. Put
$$
\zeta_3:=\zeta_2/r+\displaystyle {\frac {\varepsilon} {r\beta}},
$$
where $\beta :=\varphi'_2(\zeta_2) \not=0$. By using the hypothesis
$\varphi'_1(\zeta_2)=0,$ we have
$$
\varphi^r(\zeta_3)=(b+\alpha\varepsilon^2+O(\varepsilon^3),
\varepsilon+O(\varepsilon^2)),
\text{ where }
\alpha = \displaystyle {\frac
{\varphi_1''(\zeta_2)}{2\beta^2}}.
$$
In other words,
$\varphi^r(\zeta_3)=(b,\varepsilon)+E(\varepsilon),$ where
$E(\varepsilon) \in
\mathbb C^2$, $E(\varepsilon)=O(\varepsilon^2)$.
%  Note that we always have the estimate
%$|E_j(\varepsilon)| \le C_j|\varepsilon|^2, C_j>0, j=1,2.$
We now
define  $B:\mathbb D \to \mathbb D$ by
$$B(\zeta)=\zeta\prod^2_{j=1}\bigg ( \displaystyle{\frac
{\zeta_j/r-\zeta}{1-\zeta \overline {\zeta_j/r}}}\bigg ),$$
and $f:\mathbb D \to \mathbb C^2$ by
$$
f(\zeta)=\displaystyle{\frac
{-B(\zeta)}{B(\zeta_3)}}E(\varepsilon).
$$
Clearly, the above definitions  imply that $|B(\zeta_3)|\ge
C_0|\varepsilon|$ ($C_0$ is a constant not depending on $\varepsilon$).
Thus
$$|f_j(\zeta)|\le \displaystyle {\frac
{C_j|\varepsilon|^2}{C_0|\varepsilon|}}=:M_j|\varepsilon|, \forall
\zeta \in \mathbb D, j=1,2.$$
Put
$$\widetilde{\varphi}=\varphi^r+f.$$
When $|\varepsilon|$ is small enough, by taking 
$r=1-\displaystyle{\frac
{\max_{j=1,2}M_j}{\text{min}_{j=1,2}s_j}}|\varepsilon|$ and using the 
above
estimates we have
$$
\|\widetilde{\varphi}_j\|_{\infty}\le
\|\varphi^r_j\|_{\infty}+\|f_j\|_{\infty} \le
1-s_j(1-r)+M_j|\varepsilon| \le 1,\; j=1,2.
$$
This means that $\widetilde{\varphi}\in \mathcal{O}(\mathbb D,\mathbb
D^2).$ On the other hand, it is easy to check that
$$
\widetilde{\varphi}(0)=z, \widetilde{\varphi}(\zeta_1/r)=(a,0),
\widetilde{\varphi}(\zeta_2/r)=(b,0),
\widetilde{\varphi}(\zeta_3)=(b,\varepsilon)
$$
and
$\log|\zeta_1/r|+\log|\zeta_2/r|+\log|\zeta_3|$ tends to
$\log|\zeta_1|+2\log|\zeta_2|$ as $\varepsilon \to 0$, and thus 
$r\to1$.
It follows that
$$
\underset {\varepsilon \to 0}{\limsup}\;L_{a0bV}^\varepsilon \le
\log|\zeta_1|+2\log|\zeta_2|.
$$

In both cases, we proved that $\underset {\varepsilon \to
0}{\limsup}\;L_{a0bV}^\varepsilon \le \log|\zeta_1|+2\log|\zeta_2|.$ By
taking the infimum over all analytic discs $\varphi$ verifying
(\ref{candidate}), it
follows that
$\underset {\varepsilon \to 0}{\limsup}\;L_{a0bV}^\varepsilon \le
L_{a0bV}$. The inequality  $\underset{\varepsilon \to 0}{\text 
{limsup}}\;
L_{a0bV}^\varepsilon(z) \le \widetilde {L}(z)$ is proved.
\end{proof}

\begin{lemma}
\label{liminf}
$I:=\liminf_{\varepsilon\to 0} L_{a0bV}^\varepsilon(z) \ge \tilde 
L_{a0bV}(z)$.
\end{lemma}

\begin{proof}
Take a sequence of analytic discs $\varphi^{\varepsilon} \in \mathcal
{O}(\mathbb D,\mathbb D^2)$ such that
$$
\varphi^{\varepsilon}(0)=z, \quad
\varphi^{\varepsilon}(\zeta_1^{\varepsilon})=(a,0),
\varphi^{\varepsilon}(\zeta_2^{\varepsilon})=(b,0),
\varphi^{\varepsilon}(\zeta_3^{\varepsilon})=(b,\varepsilon)
$$
for every
$\varepsilon,$ where $\sum_{j=1}^3 \log |\zeta_j^{\varepsilon}|$
converges to $I$ as $\varepsilon \to 0.$

By passing to a subsequence, we may assume that $\varphi^{\varepsilon}$
converges locally uniformly to some $\varphi \in \mathcal {O}(\mathbb
D,\mathbb D^2).$ Also (if necessary, by passing to a subsequence
again), we may assume that $\zeta^{\varepsilon}_j \to \zeta_j \in
\overline{\mathbb D}$ for each $j$ as $\varepsilon \to 0$.

Denote $K=\{k \in \{1,2,3\}: \zeta_k \in \mathbb D\}.$ It is easy to
see that $\zeta_1 \not= \zeta_2, \zeta_1 \not= \zeta_3$ and $K \not=
\emptyset$ : if every $\zeta_j$ was in $\partial \mathbb D,$
  this would imply that $I=0$.

If either $2$ or $3 \notin K $ or $2,3 \in K$ but $\zeta_2 \not=
\zeta_3,$ then $I=\sum_{k \in K}\log|\zeta_k|$. And we have $\varphi_2
\in \mathcal{O}(\mathbb D,\mathbb D)$ with $\varphi_2(0)=z_2$ and $
\varphi_2(\zeta_k)=0, k \in K$. Hence the function $\varphi_2$ must be
of the form
$$
\varphi_2(\zeta)=\prod_{k\in K} \bigg (\displaystyle{\frac
{\zeta_k-\zeta}{1-\zeta \overline {\zeta_k}}}\bigg )h(\zeta),
$$
where $h \in \mathcal{O}(\mathbb D,\overline {\mathbb D})$ and
$h(0)=\displaystyle{\frac {z_2} {\prod_{k \in K}\zeta_k}}.$ This
implies that $|z_2|\le \prod_{k \in K}|\zeta_k|.$ Hence
$L_{a0b0}(z)=\log|z_2|\le I=\sum_{k\in K}\log|\zeta_j|.$ This proves 
that
$I \ge \widetilde {L}(z).$

If $1\notin K, 2,3\in K$ and $\zeta_2=\zeta_3.$ Then, of course, we
have $\varphi'_1(\zeta_2)=0.$ Thus $I\ge L_{bV},$ and hence $I\ge
L_{a0bV}$ by Proposition \ref{pro:minusonepoint}. So, we have the
inequality $I\ge  \widetilde {L}(z).$

If $K=\{1,2,3\}$ and $\zeta_2=\zeta_3.$ Then the function $\varphi$
belong to the defining familly of the function $L_{a0bV}.$ It implies
that $I \ge L_{a0bV}(z) \ge \tilde L_{a0bV}(z)$.
\end{proof}

This completes the proof
of Theorem \ref{liml3bepsilon=}.
\end{proof}

\section{The main counterexample}
\label{mcex}

\begin{prop}\label{prop:limsupl4epsilon}
For every $z\in \mathbb D^2$ we have
$$
\underset {\varepsilon \to
0}{\limsup}\;L_{aVbV}^\varepsilon(z)\le
\tilde L_{aVbV}(z):=\min 
\{L_{a0b0}(z),L_{aVb0}(z),L_{a0bV}(z),L_{aVbV}(z)\}.
$$
\end{prop}

\begin{proof}  By the monotonicity property of
Coman's Lempert function \cite[Theorem
2.4]{WikstromAMS} and Theorem \ref{liml3bepsilon=}, we have
$$
\underset {\varepsilon \to 0}{\limsup}\;L_{aVbV}^\varepsilon(z)\le
\underset {\varepsilon \to
0}{\limsup}\; L_{aVb0}^\varepsilon(z)=\min\{L_{a0b0}(z),L_{aVb0}(z)\},
$$
and
$$
\underset {\varepsilon \to 0}{\limsup}\;L_{aVbV}^\varepsilon(z)\le
\underset {\varepsilon \to
0}{\limsup}\;L_{a0bV}^\varepsilon(z)=\min\{L_{a0b0}(z),L_{a0bV}(z)\}.
$$
It implies that $\underset {\varepsilon \to
0}{\limsup}\;L_{aVbV}^\varepsilon(z)\le
\text{min}\{L_{a0b0}(z),L_{aVb0}(z),L_{a0bV}(z)\}$, for any $z\in 
\mathbb
D^2$.
Hence the problem has been reduced the proving $\underset {\varepsilon
\to 0}{\limsup}\;L_{aVbV}^\varepsilon(z)\le L_{aVbV}(z)$.

To prove this, we will use an argument as in the proof of Lemma
\ref{limsup}.

{\it Case 1.} $\log|z_2| \le \log \bigg |
\displaystyle {\frac {a - z_1} {1 - z_1\overline{a}}}\bigg | + \log
\bigg | \displaystyle {\frac {b - z_1} {1 - z_1\overline{b}}}\bigg |.$

A well-known special case of Theorem \ref{thm:G=LinthecaseNpoles} then
implies that
$$
L_{a0b0}(z)=\log \bigg | \displaystyle {\frac {a - z_1} {1 -
z_1\overline{a}}}\bigg | + \log \bigg | \displaystyle {\frac {b - z_1}
{1 - z_1\overline{b}}}\bigg |
=
L_{aV}(z)+L_{bV}(z)
\le
L_{aVbV}(z)
$$
Since, again by monotonicity,
  $L_{aVbV}^\varepsilon(z)\le L_{a0b0}(z)$ for every $z\in \mathbb
D^2$, we have
$$
\underset {\varepsilon \to
0}{\limsup}L_{aVbV}^\varepsilon(z)\le L_{aVbV}(z).
$$

{\it Case 2.}
$\log|z_2| > \log \bigg |
\displaystyle {\frac {a - z_1} {1 - z_1\overline{a}}}\bigg | + \log
\bigg | \displaystyle {\frac {b - z_1} {1 - z_1\overline{b}}}\bigg |$.

Then $L_{a0b0}(z)=\log|z_2|.$
Let $ \varphi \in \mathcal {O}(\mathbb D,\mathbb D^2)$ be an analytic
disc, and $\zeta_j \in \mathbb D$, $j=1,2$ such that
\begin{equation}
\label{cand2}
\varphi(0)=z,
\varphi(\zeta_1)=(a,0), \varphi(\zeta_2)=(b,0),
\varphi_1'(\zeta_1)=\varphi'_1(\zeta_2)=0.
\end{equation}

We consider the following three cases.

{\it Case 2.1} $\varphi_2'(\zeta_j)\not= 0$, $j=1,2$.

Fix $\varepsilon
\in \mathbb C$ with $|\varepsilon|$ small enough. We will show that
there exist $\widetilde{\varphi}\in \mathcal{O}(\mathbb D,\mathbb 
D^2)$,
$\widetilde{\zeta_j} \in \mathbb D$, $1 \le j \le 4$ such that
$$
\widetilde{\varphi}(\widetilde{\zeta_1})=(a,0), \;
\widetilde{\varphi}(\widetilde{\zeta_2})=(b,0), \;
\widetilde{\varphi}(\widetilde{\zeta_3})=(b,\varepsilon), \;
\widetilde{\varphi}(\widetilde{\zeta_4})=(a,\varepsilon)
$$
and
$$
\sum_{j=1}^4\log|\widetilde{\zeta_j}| \to
2\log|\zeta_1|+2\log|\zeta_2|  \; \text{as}\; \varepsilon \to 0.
$$

Let $1/2 \le r<1$ to be chosen later. As in the proof of Lemma
\ref{limsup}, we define $\varphi^r:\mathbb D \to \mathbb D^2$
by $\varphi^r(\zeta)=\varphi(r\zeta)$. From (\ref{cand2}) we have
$$
\varphi^r(0)=z, \varphi^r(\zeta_1/r)=(a,0),
\varphi^r(\zeta_2/r)=(b,0),
$$
and, as in the proof of Lemma \ref{limsup},
$|\varphi_j^r(\zeta)| \le 1-s_j
(1-r)$, for any $ \zeta \in \mathbb D, j=1,2,$ where $s_j$ depends
only on $z$.

We set
\begin{eqnarray*}
\widetilde{\zeta_0}&:=&0,\qquad \widetilde{\zeta_1}:=\zeta_1/r, \qquad
\widetilde{\zeta_2}:=\zeta_2/r,\\
\widetilde{\zeta_3} &:=&
\frac{\zeta_2}r+\displaystyle {\frac {\varepsilon} {r\beta_2}},\;
\text{where}\; \beta_2=\varphi'_2(\zeta_2) \not=0, \\
\widetilde{\zeta_4}&:=&\frac{\zeta_1}r+\displaystyle {\frac {\varepsilon}
{r\beta_1}},\;\text{where} \;\beta_1=\varphi'_2(\zeta_1) \not=0.
\end{eqnarray*}

Using the hypothesis  $\varphi_1'(\zeta_j)=0, j=1,2$ we have
$$
\varphi^r(\widetilde{\zeta_3})=(b,\varepsilon)+E^2(\varepsilon),
\;
\text{where}\;E^2(\varepsilon)
=(\alpha_2\varepsilon^2+O(\varepsilon^3),O(\varepsilon^2)) \in
\mathbb C^2, \alpha_2 = \displaystyle {\frac
{\varphi_1''(\zeta_2)}{2\beta_2^2}},
$$
$$
\varphi^r(\widetilde{\zeta_4})=(a,\varepsilon)+E^1(\varepsilon),
\;
\text{where}\;E^1(\varepsilon)
=(\alpha_1\varepsilon^2+O(\varepsilon^3),O(\varepsilon^2))\in
\mathbb C^2, \alpha_1 = \displaystyle {\frac
{\varphi_1''(\zeta_1)}{2\beta_1^2}},$$
and of course $|E^k_j(\varepsilon)| \le C^k_j|\varepsilon|^2, C^k_j>0,
k,j=1,2.$ We now define  $f:\mathbb D \to \mathbb C^2$ by
$$
f(\zeta)
=
-E^1(\varepsilon)\prod^3_{j=0}\bigg ( \displaystyle{\frac
{\zeta-\widetilde{\zeta_j}}{\widetilde{\zeta_4}-\widetilde{\zeta_j}}}\bigg
)
-E^2(\varepsilon)\prod^4_{j=0,j\not=3}\bigg ( \displaystyle{\frac
{\zeta-\widetilde{\zeta_j}}{\widetilde{\zeta_3}-\widetilde{\zeta_j}}}\bigg
).
$$
It is easy to see that $|f_j(\zeta)|\le
M_j|\varepsilon|, \forall \zeta \in \mathbb D, j=1,2,$ where $M_j$ are
positive constants which do not depend on $r$.

Put $ \widetilde{\varphi}(\zeta):=\varphi^r(\zeta)+f(\zeta)$
and $r=1-\displaystyle{\frac
{\max_{j=1,2}M_j}{\text{min}_{j=1,2}s_j}}|\varepsilon|$. Then
$$
\widetilde{\varphi}(0)=z, \;
\widetilde{\varphi}(\widetilde{\zeta_1})=(a,0), \;
\widetilde{\varphi}(\widetilde{\zeta_2})=(b,0), \;
\widetilde{\varphi}(\widetilde{\zeta_3})=(b,\varepsilon), \;
\widetilde{\varphi}(\widetilde{\zeta_4})=(a,\varepsilon) ;
$$
$$\mbox{and }
\|\widetilde{\varphi}_j\|_{\infty}\le
\|\varphi^r_j\|_{\infty}+\|f_j\|_{\infty} \le
1-s_j(1-r)+M_j|\varepsilon| \le 1,\; j=1,2
$$
when $|\varepsilon|$ is small enough.
This means
that $\widetilde{\varphi}\in \mathcal{O}(\mathbb D,\mathbb D^2)$.
Finally
$\sum_{j=1}^4\log|\widetilde{\zeta_j}|$ tends to
$2\log|\zeta_1|+2\log|\zeta_2|$ as $\varepsilon \to 0$. This proves
that $\underset {\varepsilon \to 
0}{\limsup}\;L_{aVbV}^\varepsilon(z)\le
2\log|\zeta_1|+2\log|\zeta_2|$.

{\it Case 2.2} One of the $\zeta_j, j=1,2,$ is a zero of the
function $\varphi'_2$.

Without loss of generality we may assume that
$\varphi_2'(\zeta_1)\not= 0, \varphi_2'(\zeta_2)=0.$ Let $0<r<1$. By
using again the function $\varphi^r$ as in the previous case, we 
already have
  $|\varphi_j^r(\zeta)| \le 1-s_j (1-r)$, for any $\zeta \in \mathbb
D$, $j=1,2$. Define $\psi: \mathbb D \to \mathbb C^2$ by
$$
\psi(\zeta):=\varphi^r(\zeta)+\bigg (0,C \cdot \zeta \cdot \bigg
(\displaystyle{\frac
{\zeta_1/r-\zeta}{1-\overline{\zeta_1/r}\zeta}}\bigg)^2\bigg
(\displaystyle{\frac
{\zeta_2/r-\zeta}{1-\overline{\zeta_2/r}\zeta}}\bigg)\bigg),
$$
where $C$ is a positive constant small enough so that $\psi \in
\mathcal {O}(\mathbb D, \mathbb D^2)$.
Clearly, $\psi$ and $\zeta_1/r$, $\zeta_2/r$ belong to the defining 
family
of the function $L_{aVbV}$ and $\psi'_2(\zeta_j)\not=0, j=1,2$. By using
again the proof of Case 2.1 with $\psi$ instead of $\varphi$ and
$\zeta_j/r$ instead of $\zeta_j,$ we have the inequality $\underset
{\varepsilon \to 0}{\limsup}\;L_{aVbV}^\varepsilon(z)\le
2\log|\zeta_1|+2\log|\zeta_2|.$

{\it Case 2.3} $\varphi_2'(\zeta_j)=0, j=1,2.$

Using an argument as
in the case 2.2 and considering the function
$$
\psi(\zeta):=\varphi^r(\zeta)+\bigg (0,C \cdot \zeta \cdot \bigg
(\displaystyle{\frac
{\zeta_1/r-\zeta}{1-\overline{\zeta_1/r}\zeta}}\bigg) \bigg
(\displaystyle{\frac
{\zeta_2/r-\zeta}{1-\overline{\zeta_2/r}\zeta}}\bigg) \bigg),
$$
we have the equality $\underset {\varepsilon \to
0}{\limsup}\;L_{aVbV}^\varepsilon(z)\le 2\log|\zeta_1|+2\log|\zeta_2|.$

Thus we always do have
$$
\underset {\varepsilon \to 0}{\limsup}\;L_{aVbV}^\varepsilon(z)\le
2\log|\zeta_1|+2\log|\zeta_2|.
$$
By taking the infimum over all  discs $\varphi$ satisfying 
(\ref{cand2}),
it follows
that
$$\underset {\varepsilon \to 0}{\limsup}L_{aVbV}^\varepsilon(z)\le
L_{aVbV},$$
and the proof is finished.
\end{proof}

\begin{thm}
\label{NoComan}
Coman's question admits a negative answer in the bidisc for equal 
weights. More
precisely, consider
$$
S_{aVbV}^\varepsilon:= \{(a,0); (b,0), (b,\varepsilon),
(a,\varepsilon)\} \mbox{ with }b =-a,
$$
where $\varepsilon \in
\mathbb C$ and the weights are all equal to $1$, and $z=(0,\gamma)$ 
with
$|a|^{3/2}<|\gamma|<|a|$. Then, 
$\underset {\varepsilon \to 0}{\liminf}
L_{aVbV}^\varepsilon > G_{aVbV}(z)$ and therefore, 
for $|\varepsilon|$ small enough,
$$
G_{aVbV}^\varepsilon(z) < L_{aVbV}^\varepsilon(z) .
$$
\end{thm}

\begin{proof}

Using the result of Edigarian about the product property of the Green
function,
  \cite{Edigarian}, we have
$$
G_{aVbV}^\varepsilon=\max\bigg
\{\log|a|+\log|b|;\log|\gamma|+\log\bigg|\displaystyle{\frac
{\varepsilon-\gamma}{1-\overline{\varepsilon}\gamma}}\bigg|\bigg \}.
$$
Thus
$$
G_{aVbV}(z)=\underset {\varepsilon \to
0}{\lim}G_{aVbV}^\varepsilon=\log|a|+\log|b|=\log|a|^2.
$$
By Propositions \ref{prop:G2ab<L2ab} and \ref{prop:G2ab<L2b}, and since
$L_{a0b0}(z)=\log|\gamma|>\log|a|^2 = G_{aVbV}(z),$ we have
\begin{equation}
\label{GVV<LVV}
\underset {\varepsilon \to
0}{\lim}G_{aVbV}^\varepsilon=G_{aVbV}(z)<\tilde L_{aVbV}(z):=\text
{min}\{L_{a0b0}(z), L_{aVb0}(z),L_{a0bV}(z),L_{aVbV}(z)\}.
\end{equation}

We consider $I:=\underset {\varepsilon \to 0}{\liminf}
L_{aVbV}^\varepsilon$. We want to prove that $I > 
G_{aVbV}(z)$.
As in the proof of Lemma \ref{liminf}
take for each $\varepsilon$ an analytic disc
  $\varphi^{\varepsilon} \in \mathcal
{O}(\mathbb D,\mathbb D^2)$ such that
$$
\varphi^{\varepsilon}(0)=z,
\varphi^{\varepsilon}(\zeta_1^{\varepsilon})=(a,0),
\varphi^{\varepsilon}(\zeta_2^{\varepsilon})=(b,0),
\varphi^{\varepsilon}(\zeta_3^{\varepsilon})=(b,\varepsilon),
\varphi^{\varepsilon}(\zeta_4^{\varepsilon})=(a,\varepsilon)
$$
and such that $\sum_{j=1}^4 \log |\zeta_j^{\varepsilon}|$
converges to $I$ as $\varepsilon \to 0$.

By passing to a subsequence, we may assume that $\varphi^{\varepsilon}$
converges locally uniformly to some $\varphi \in \mathcal {O}(\mathbb
D,\mathbb D^2).$ Also (if necessary, by passing to a subsequence
again), we may assume that $\zeta^{\varepsilon}_j \to \zeta_j \in
\overline{\mathbb D},$ for each $j,$ as $\varepsilon \to 0.$

Denote $K=\{k \in \{1,2,3,4\}: \zeta_k \in \mathbb D\}.$ It is easy to
see that $\mathbb D \cap \{ \zeta_1, \zeta_4\} \cap \{ 
\zeta_2,\zeta_3\}
= \emptyset$.

If $K=\emptyset$ then $I=0,$ and hence we have
$I\ge \tilde L_{aVbV}(z)>G_{aVbV}(z)$, by (\ref{GVV<LVV}).
So now we only consider the cases where
$K\not=\emptyset.$

If $\zeta_j \not= \zeta_k$, $\forall j \neq k \in K$, then $I=\sum_{k 
\in
K}\log|\zeta_k|$, $\varphi_2 \in \mathcal{O}(\mathbb D,\mathbb D)$,
 $\varphi_2(0)=\gamma$ and $\varphi_2(\zeta_k)=0, k \in K.$ It
implies that
$$
\varphi_2(\zeta)=\prod_{k\in K} \bigg (\displaystyle{\frac
{\zeta_k-\zeta}{1-\zeta \overline {\zeta_k}}}\bigg )h(\zeta),
$$
where $h \in \mathcal{O}(\mathbb D,\overline {\mathbb D})$ and
$h(0)=\displaystyle{\frac {\gamma} {\prod_{k \in K}\zeta_k}}.$ Thus we
have
$$
L_{a0b0}(z)= \log|\gamma| \le \sum_{k\in K}\log|\zeta_k|=I,
$$
and
hence, $I \ge \tilde L_{aVbV}(z)>G_{aVbV}(z).$

If $K=\{2,3\}$ and $\zeta_2=\zeta_3,$  then, since
$\zeta^{\varepsilon}_2 \to \zeta_2$,  $\zeta^{\varepsilon}_3 \to 
\zeta_2$
and $|\zeta^{\varepsilon}_3-\zeta^{\varepsilon}_2| \ge |\varepsilon|$,
$$
\varphi'_1(\zeta_2)=
\lim_{\varepsilon \to 0}
\frac{0}{\zeta^{\varepsilon}_3-\zeta^{\varepsilon}_2} =
0.
$$
Thus $I\ge L_{bV}(z)\ge L_{a0bV}(z)$ by Proposition
\ref{pro:minusonepoint}. So that $I \ge \tilde L_{aVbV}(z)> 
G_{aVbV}(z).$

If $K=\{1,4\}$ and $\zeta_1=\zeta_4,$  then $\varphi'_1(\zeta_1)=0.$
Thus $I\ge L_{aV}(z)\ge L_{aVb0}(z)$ by Proposition
\ref{pro:minusonepoint}. So that $I \ge \tilde L_{aVbV}(z)> 
G_{aVbV}(z)$.

If $K=\{1,2,3\}, \zeta_2=\zeta_3$,
then $\varphi'_1(\zeta_2)=0.$ Thus $I= \log|\zeta_1| + 2 \log |\zeta_2| \ge 
L_{a0bV}(z)\ge
\tilde L_{aVbV}(z)> G_{aVbV}(z).$ The same reasoning obtains if 
or  $K=\{4,2,3\}$,
$\zeta_2=\zeta_3$.

Similarly,
if either $K=\{1,2,4\}, \zeta_1=\zeta_4$ or  $K=\{1,3,4\},
\zeta_1=\zeta_4,$ then $\varphi'_1(\zeta_1)=0.$ This implies that $I\ge
L_{aVb0}(z)\ge \tilde L_{aVbV}(z)> G_{aVbV}(z).$

If $K=\{1,2,3,4\}$ and $\zeta_1=\zeta_4, \zeta_2=\zeta_3,$ then $
\varphi'_1(\zeta_1)=\varphi'_1(\zeta_2)=0.$ It implies that  
$I= 2\log|\zeta_1| + 2 \log |\zeta_2|\ge
L_{aVbV}(z)\ge \tilde L_{aVbV}(z)> G_{aVbV}(z).$

Suppose now that
 $K=\{1,2,3,4\}$ and $\zeta_1\not=\zeta_4$, $\zeta_2=\zeta_3.$ 
This is the final and most delicate case ; the proof of Theorem 
\ref{limell4<ltilde} 
below suggests that it may occur for some values of $\gamma$. Both previous 
types of argument now break down, because we only get 
$$
I < \min(\log|\zeta_1|,\log|\zeta_4|) + 2 \log |\zeta_2|\ge L_{a0bV}(z) ;
$$
or, from the fact that 
$\varphi_2(\zeta_1)=\varphi_2(\zeta_4)=\varphi_2(\zeta_2)= 0$ and 
$\varphi_2(0)= \gamma$, 
$$
I < \log|\zeta_1| + \log|\zeta_4| + \log |\zeta_2| \ge \log| \gamma| 
\ge L_{a0b0}(z) .
$$

By using
a rotation in the first coordinate we can assume that $a>0.$ We 
will prove that $I>G_{aVbV}(z).$ If not, we would have
\begin{equation}
\label{assume}
\log|\zeta_1|+\log|\zeta_4|+2\log|\zeta_2|=I=G_{aVbV}(z)=2\log a.
\end{equation}
Then the function $\varphi_1$  has the following properties:
\begin{equation}
\label{require1}
\varphi_1(0)=0;  \;
\varphi_1(\zeta_1)=\varphi_1(\zeta_4)=a; \;
\varphi_1(\zeta_2)=-a;  \;
\varphi'_1(\zeta_2)=0.
\end{equation}

Setting $f:=\phi_{-a}\circ\varphi_1\circ\phi_{\zeta_2}$, with 
$\phi_\xi$
defined as in (\ref{Mobmap}),
we have $f(0)=0$, $f'(0)=0$ and $f(\zeta_2)=-a$. The Schwarz
Lemma shows that $|\zeta_2|^2\ge a,$ and hence
\begin{equation}
\label{zetatwo}
2\log|\zeta_2| \ge \log\, a
\end{equation}
Setting $g:=\phi_{a} \circ \varphi_1$, we have
$g(\zeta_1)=g(\zeta_4)=0$ and $g(0)=a$. Thus the function $g$ must have
the following form
$$
g(\zeta)=\phi_{\zeta_1}(\zeta)\phi_{\zeta_4}(\zeta)h_1(\zeta),
\forall \zeta \in \mathbb D,
  \text{ where } h_1 \in
\mathcal{O}(\mathbb D,\overline {\mathbb D})
\text{ and }
h_1(0)=\displaystyle{\frac {a} {\zeta_1\zeta_4}}, \mbox{ hence }
$$
\begin{equation}
\label{zetaonefour}
\log|\zeta_1|+\log|\zeta_4| \ge \log a .
\end{equation}
 From (\ref{zetatwo}) and (\ref{zetaonefour}) we have
$$
I=\log|\zeta_1|+\log|\zeta_4|+2\log|\zeta_2| \ge 2\log\; a.
$$
The assumption (\ref{assume}) implies that all the inequalities in
(\ref{zetatwo}) and (\ref{zetaonefour}) become  equalities.
Now, since $\varphi_2(0)=\gamma$ and
$\varphi_2(\zeta_1)=\varphi_2(\zeta_2)=\varphi_2(\zeta_4)=0$,
  $$
  \varphi_2(\zeta)=\prod_{j=1, j\not= 3}^4\bigg (\displaystyle {\frac
{\zeta_j-\zeta}{1-\overline{\zeta_j}\zeta}}\bigg )h_2(\zeta),
\quad \text{ where } h_1 \in
\mathcal{O}(\mathbb D,\overline {\mathbb D})
\text{ and }
h(0)=\displaystyle{\frac {\gamma}{\zeta_1\zeta_2\zeta_4}}.
$$
It implies
that $|\gamma|\le |\zeta_1\zeta_2\zeta_4|=a^{3/2}.$ This contradicts
the hypothesis $|\gamma| > a^{3/2},$ and the inequality $I>G_{aVbV}(z)$
is proved.

If $K=\{1,2,3,4\}$ and $\zeta_1=\zeta_4$, $\zeta_2\not=\zeta_3$, the 
proof
is similar.
\end{proof}

We now prove that the function $\tilde L_{aVbV}$ is not, in spite of 
the
positive result obtained in Proposition \ref{prop:limsupl4epsilon}, the
limit of the functions $L_{aVbV}^\varepsilon$.

\begin{thm}\label{limell4<ltilde}
Let $S_{aVbV}^\varepsilon$ be as in Theorem \ref{NoComan}, with $b 
=-a$,
and $z=(0,\gamma)$ with
$|a|^2<|\gamma| \leq |a|^{3/2}$. Then, for $|\varepsilon|$ small 
enough,
$$
  \underset {\varepsilon \to 0}{\limsup} L_{aVbV}^\varepsilon<
\tilde L_{aVbV}.
$$
\end{thm}

\begin{proof} Without loss of generality we may assume that $a>0$. Then
we have $G_{aVbV}(z)=2\log a < \log|\gamma|=L_{a0b0}(z)$. 
By Propositions \ref{prop:G2ab<L2ab} and
\ref{prop:G2ab<L2b}, we have $G_{aVbV}(z) < \min( L_{aVbV}(z), 
L_{a0bV}(z),  L_{aVb0}(z) )$, so $G_{aVbV}(z) < \tilde L_{aVbV}(z)$. 

We now prove that $\underset {\varepsilon \to 0}{\limsup} L_{aVbV}^\varepsilon
\le G_{aVbV}(z)$
in two steps. In Step 1 we construct $\varphi \in \mathcal
{O}(\mathbb D,\mathbb D^2)$ and $\zeta_1$, $\zeta_2$, $\zeta_4 \in 
\mathbb
D$ such that
\begin{equation}
\label{requirephi}
\varphi(0)=z= (0,\gamma), \quad
\varphi (\zeta_1)=\varphi(\zeta_4)=(a,0), \quad
\varphi(\zeta_2)=(b,0), \quad
\varphi'_1(\zeta_2)=0
\end{equation}
and
\begin{equation}
\label{requirezeta}
\log|\zeta_1|+\log|\zeta_4|+2\log|\zeta_2|=2\log a.
\end{equation}
In Step 2, with the methods of the proofs of Lemma \ref{limsup} and
Proposition \ref{prop:limsupl4epsilon},
we prove that for any $\varepsilon \in \mathbb C$
with $|\varepsilon|$ small enough, there are $\widetilde{\varphi}\in
\mathcal{O}(\mathbb D,\mathbb D^2)$, $\widetilde{\zeta_j} \in \mathbb 
D$,
$1 \le j \le 4$ such that
$$
\widetilde{\varphi}(0)= z,  \quad
\widetilde{\varphi}(\widetilde{\zeta_1})=(a,0), \quad
\widetilde{\varphi}(\widetilde{\zeta_2})=(b,0), \quad
\widetilde{\varphi}(\widetilde{\zeta_3})=(b,\varepsilon), \quad
\widetilde{\varphi}(\widetilde{\zeta_4})=(a,\varepsilon)
$$
and
\begin{equation}
\label{limofsum}
\sum_{j=1}^4\log|\widetilde{\zeta_j}| \to
\log|\zeta_1|+\log|\zeta_4|+2\log|\zeta_2|=2\log \, a \; \text{as}\;
\varepsilon \to 0.
\end{equation}
This shows that $\underset {\varepsilon \to
0}{\limsup}L_{aVbV}^\varepsilon(z)\le 2\log \, a=G_{aVbV}(z).$ 

{\it Step 1.} Let $\varphi=(\varphi_1,\varphi_2): \mathbb D \to \mathbb
D^2$ satisfying (\ref{requirephi}), then $\phi_1$ satisfies
(\ref{require1}) and as in the proof of Theorem \ref{NoComan}, we see
that $\zeta_2$ verifies (\ref{zetatwo}) and $\zeta_1$ and $\zeta_4$
verify (\ref{zetaonefour}). From this and (\ref{requirezeta}), we see 
that
all inequalities in (\ref{zetatwo}) and (\ref{zetaonefour}) become
equalities.

By choosing
$\zeta_2=\sqrt{a}$,  we have to have
$$
f(\zeta) :=
\phi_{-a}\circ \varphi_1\circ \phi_{\zeta_2}(\zeta)=- \zeta^2,
\quad \forall \zeta \in
\mathbb D.
$$
Hence
$$
\varphi_1(\zeta)=\phi_{-a}\bigg (-\phi^2_{\zeta_2}(\zeta)\bigg).
$$
Using the conditions
$f(\phi_{\zeta_2}(\zeta_1))=f(\phi_{\zeta_2}(\zeta_4))=\phi_{-a}(a)$, 
we
get that
$$
\zeta_1=\phi_{\sqrt{a}}(\xi); \zeta_4=\phi_{\sqrt{a}}(-\xi),
\mbox{ where } \xi := \sqrt{\displaystyle {\frac {2a}{1+a^2}}}.
$$
To get a function $\varphi_2$ such that
$\varphi_2(\zeta_j)=0$, $j=1,2,4$ and $\varphi_2(0)=\gamma$, noting
that $|\zeta_1\zeta_2\zeta_4|= |a|^{3/2} >|\gamma|>0$, set
$$
\varphi_2(\zeta)=\frac {\gamma}{\zeta_1\zeta_2\zeta_4}
\prod^4_{j=1, j\not=3}\phi_{\zeta_j}(\zeta), \quad
\forall \zeta \in \mathbb D .
$$
Then $\varphi_2'(\zeta_2)\not=0$. The
analytic disc $\varphi =(\varphi_1,\varphi_2)$ is well defined and
satisfies all the required properties.

{\it Step 2.}
Let $0<r<1$ be chosen later. As in the proof of Theorem
\ref{liml3bepsilon=} and Proposition \ref{prop:limsupl4epsilon}, we
define $\varphi^r:\mathbb D \to \mathbb D^2$ by
$\varphi^r(\zeta)=\varphi(r\zeta)$. It is
easy to see that
$$
\varphi^r(0)=z, \;
\varphi^r(\zeta_1/r)=(a,0), \;
\varphi^r(\zeta_2/r)=(-a,0),  \;
\left(\varphi^r_1\right)'(\zeta_2/r)= 0, \;
\varphi^r(\zeta_4/r)=(a,0)
$$
and
$|\varphi_j^r(\zeta)| \le 1-s_j (1-r),\forall \zeta \in \mathbb D,
j=1,2,$ where the $s_j$ depend only on $z$.

Setting
  $$
\zeta_3:=\zeta_2/r+\displaystyle {\frac {\varepsilon} {r\beta}},
\text{ where } \beta=\varphi'_2(\zeta_2) \not=0 ,
$$
by using the hypotheses  $\varphi_1'(\zeta_j)=0$, $j=1,2$ we have
$$
\varphi^r(\zeta_3)=(b,\varepsilon)+E(\varepsilon),\text{ where }
E(\varepsilon)=O(\varepsilon^2)
$$
As in the proof of Proposition \ref{prop:limsupl4epsilon}
we can define the function
$f:\mathbb D \to \mathbb C^2$ which interpolates at the points
  $\{0, \zeta_1/r, \zeta_2/r, \zeta_3, \zeta_4/r \}$ in the unit disc 
$\mathbb
d$ respectively the values
$$
\left\{ (0,0), (0,0), (0,0), -E(\varepsilon),
(0,\varepsilon) \right\}
$$ in $\mathbb C^2$ and
$|f_j(\zeta)|\le
M_j|\varepsilon|$, for all $\zeta \in \mathbb D$, $j=1,2,$ where the 
$M_j$ are
positive constants.

Setting
$\widetilde{\varphi}(\zeta):=\varphi^r(\zeta)+f(\zeta),$
and $\widetilde{\zeta}_j := \zeta_j/r$, $j=1,2,4$, 
$\widetilde{\zeta}_3 := \zeta_3$, and $|\varepsilon|$ small enough
with respect to $1-r$,
all required properties are verified. Since $1-r$ is arbitrarily small,
we do have (\ref{limofsum}).
\end{proof}

\noindent {\bf Acknowledgements.} The results of this paper will be part of 
the second-named author's Ph. D. dissertation, and were obtained in part 
during a stay at the Paul
Sabatier university. He would like to thank Professor Do Duc Thai
for  stimulating discussions regarding this paper, and 
 Professors Nguyen Thanh Van and Fr\'ed\'edric Pham 
for the invitation and financial support, and  the Emile Picard
laboratory for hospitality.  The first-named author thanks St\'ephanie 
Nivoche and Evgueny Poletsky for interesting discussions concerning the 
topic, and also the latter for his hospitality in Syracuse.
\vspace{0.5cm}

% ----------------------------------------------------------------
\bibliographystyle{amsplain}

\vskip1cm
%\vfill\eject

Pascal J. Thomas

Laboratoire Emile Picard, UMR CNRS 5580

Universit\'e Paul Sabatier

118 Route de Narbonne

F-31062 TOULOUSE CEDEX

France

pthomas@cict.fr

\vskip.5cm

Nguyen Van Trao

Department of Mathematics

Dai Hoc Su Pham 1 (Pedagogical Institute of Hanoi)

Cau Giay, Tu Liem

Ha Noi

Viet Nam

ngvtrao@yahoo.com

\end{document}